\documentstyle[mssymb]{article}

\catcode`\@=11

\font\teneuf=eufm10
\font\seveneuf=eufm7
\font\fiveeuf=eufm5
\newfam\euffam
\textfont\euffam=\teneuf
\scriptfont\euffam=\seveneuf
\scriptscriptfont\euffam=\fiveeuf
\def\next{\Err@{Use \string\frak\space only in math mode}}
\def\frak{\ifmmode\let\next=\frak@\else
 \def\next{\Err@{Use \string\frak\space only in math mode}}\fi\next}
\def\frak@#1{{\frak@@{#1}}}
\def\frak@@#1{\fam\euffam#1}

\catcode`\@=12

\newpage
\newcommand{\labels}[2]{}
\newcommand{\aftersection}{\setcounter{equation}{0}\setcounter{paragraph}{0}}

\newcommand{\bse}[2]{\section{#1}\label{se:#2}\aftersection\labels{
}{\mnote{#2}}}
\newcommand{\bsu}[2]{\paragraph{#1.}\label{su:#2}\labels{}{\mnote{#2}}}

\newtheorem{proposition}{\sc Proposition}[section]
\newtheorem{theorem}[proposition]{\sc Theorem}
\newtheorem{lemma}[proposition]{\sc Lemma}
\newtheorem{remark}[proposition]{\sc Remark}
\newtheorem{definition}[proposition]{\sc Definition}
\newtheorem{corollary}[proposition]{\sc Corollary}

\newenvironment{proof}{\begin{list}{}{\leftmargin=0pt\rightmargin=0pt}
\item {\em Proof.\ }}{\rule{6pt}{6pt}\end{list}\vspace{4pt plus 2pt minus
1pt}}

\newcommand{\bth}[1]{\begin{theorem}\label{th:#1}\labels{}{\mnote{#1}}}
\newcommand{\bthn}[2]{\begin{theorem}[#1]\label{th:#2}\labels{}{\mnote{#2}}}
\renewcommand{\eth}{\end{theorem}}
\newcommand{\bpr}[1]{\begin{proposition}\label{pr:#1}\labels{}{\mnote{#1}}}
\newcommand{\epr}{\end{proposition}}
\newcommand{\bco}[1]{\begin{corollary}\label{co:#1}\labels{}{\mnote{#1}}}
\newcommand{\eco}{\end{corollary}}
\newcommand{\ble}[1]{\begin{lemma}\label{le:#1}\labels{}{\mnote{#1}}}
\newcommand{\ele}{\end{lemma}}
\newcommand{\bre}[1]{\begin{remark}\label{re:#1}\labels{}{\mnote{#1}}\rm}
\newcommand{\ere}{\end{remark}}
\newcommand{\bde}[1]{\begin{definition}\label{de:#1}\labels{}{\mnote{#1}}\rm}
\newcommand{\ede}{\end{definition}}
\newcommand{\de}[1]{{\sc #1}}
\newcommand{\bpf}{\begin{proof}}
\newcommand{\epf}{\end{proof}}

\newcommand{\refpr}[1]{Proposition \ref{pr:#1}}
\newcommand{\refre}[1]{Remark \ref{re:#1}}
\newcommand{\refth}[1]{Theorem \ref{th:#1}}
\newcommand{\refco}[1]{Corollary \ref{co:#1}}
\newcommand{\refle}[1]{Lemma \ref{le:#1}}

\newcommand{\refeq}[1]{(\ref{eq:#1})}
\newcommand{\refse}[1]{\S\ref{se:#1}\,\,}
\newcommand{\refsu}[1]{\S\ref{su:#1}\,\,}
\newcommand{\refpa}[1]{part \ref{#1}}

\newcommand{\beql}[1]{\labels{}{\mnote{#1}}\begin{eqnarray}\label{eq:#1}}
\newcommand{\eeql}{\end{eqnarray}}
\newcommand{\beq}{\begin{eqnarray*}}
\newcommand{\eeq}{\end{eqnarray*}}

\newcommand{\al}{\alpha}

\newcommand{\la}{\lambda}

\newcommand{\G}{\Gamma}

\newcommand{\La}{\Lambda}

\newcommand{\C}{{\Bbb C}}

\renewcommand{\P}{{\Bbb P}}

\newcommand{\R}{{\Bbb R}}
\newcommand{\Z}{{\Bbb Z}}

\newcommand{\cd}{{\cal D}}

\newcommand{\cg}{{\cal G}}
\newcommand{\ch}{{\cal H}}

\newcommand{\cm}{{\cal M}}
\newcommand{\cn}{{\cal N}}
\newcommand{\co}{{\cal O}}

\newcommand{\cs}{{\cal S}}
\newcommand{\ct}{{\cal T}}

\newcommand{\hymh}{Hermitian-Yang-Mills-Higgs}
\newcommand{\hym}{Hermitian-Yang-Mills}
\newcommand{\ymh}{Yang-Mills-Higgs}
\newcommand{\ie}{{\rm i.\,e.\ }}
\newcommand{\eg}{{\rm e.\,g.\ }}
\newcommand{\vb}{$V$-bundle}

\newcommand{\dbar}{\o\partial}

\renewcommand{\deg}{\mbox{deg\,}}

\newcommand{\tr}{\mbox{tr\,}}
\newcommand{\Herm}{{\rm Herm\,}}
\newcommand{\End}{{\rm End\,}}

\newcommand{\Aut}{{\rm Aut\,}}

\newcommand{\Hom}{{\rm Hom\,}}

\renewcommand{\and}{\quad\mbox{ and }\quad}

\renewcommand{\o}[1]{\overline{#1}}
\newcommand{\oo}[1]{\overline{#1}'}

\renewcommand{\b}[1]{\breve{#1}}

\newcommand{\ci}{{\rm i}}
\newcommand{\isum}{\sum_{i=1}^n}
\newcommand{\ilist}[1]{#1_1,\dots,#1_n}

\newcommand{\supsetne}{\setbox0=\hbox{$\ne$}%
                     \setbox1=\hbox{$ni$}%
                     \supset\kern-\wd1
                     \lower1.2ex\box0}
\newcommand{\surjarrow}{\setbox0=\hbox{$\rightarrow$}%
                     \rightarrow\kern-\wd0
                     \longrightarrow}

\newcommand{\bprn}[2]{\begin{proposition}[#1]\label{pr:#2}}
\newcommand{\wo}[1]{\widetilde{#1}}
\newcommand{\lift}{\widehat}

\begin{document}

\title{{\bf Orbifold Riemann Surfaces and the \ymh\ Equations}} \author{{\bf
Ben
Nasatyr\thanks{Current address: Department of Mathematical Sciences,
University of Aberdeen,
Edward Wright Building, Dunbar Street, Old Aberdeen, AB9 2TY.}\,
and Brian
Steer,}\\ Peterhouse, Cambridge and Hertford College, Oxford} \date{Preprint,
August 6,
1993,\\
revised, January 5 and April 26, 1995}

\maketitle %\tableofcontents

\pagenumbering{arabic} %\newpage

\setcounter{section}{-1}

\bse{Introduction}{int}

In this paper we study the $U(2)$ \ymh\ equations on orbifold Riemann surfaces.
Among other aspects, we discuss existence theorems for solutions of the \ymh\
equations, the analytic construction of the moduli space of such solutions, the
connectivity and topology of this space, its holomorphic symplectic structure
and its reinterpretations as a space of orbifold Higgs bundles or
$SL_2(\C)$-representations of (a central extension of) the orbifold fundamental
group.  We follow Hitchin's original paper for (ordinary) Riemann surfaces
\cite{hi87} quite closely but there are many novelties in the orbifold
situation.  (There is some overlap with a recent paper of Boden and
Yokogawa \cite{by}.)

It may help to mention here a few of our motivations.
\begin{enumerate}
\item In studying the orbifold moduli space, we are also studying the parabolic
moduli
space (see \refsu{parhig}, also \cite{si90,by,ns93}).
\item The moduli space provides interesting examples of non-compact
hyper-K\"ahler
manifolds in all dimensions divisible by 4.
\item As a special case of the existence theorem for solutions of the \ymh\
equations
we
get the existence of metrics with conical singularities and constant sectional
curvature
on `marked' Riemann surfaces (see \refco{negative curvature}, \refth{conical}
and
compare \cite{ht92}).
\item The orbifold fundamental groups we study are Fuchsian groups and their
central
extensions: these include the fundamental groups of elliptic surfaces and of
Seifert
manifolds. We obtain results on
varieties of $SL_2(\C)$- and $SL_2(\R)$-representations of such
groups (see \refse{rep} and compare \eg
\cite{jn85}). In particular, we prove that Teichm\"uller space for a Fuchsian
group or, equivalently, for a `marked' Riemann surface is homeomorphic to
a ball  (\refth{ball}).
\item Moduli of parabolic Higgs bundles and of marked Riemann surfaces have
potential
applications in Witten's work on Chern-Simons gauge theory.
\end{enumerate}

Let $E$ be a Hermitian rank 2 $V$-bundle (\ie orbifold bundle) over an orbifold
Riemann
surface of negative Euler characteristic, equipped with a normalised volume
form,
$\Omega$. Let $A$ be a unitary connexion on $E$ and $\phi$ an
$\End(E)$-valued $(1,0)$-form. Then the \ymh\ equations are
\beq
\begin{array}{rcl} F_A + [\phi,\phi^*] &=& - \pi \ci\, c_1(E)\Omega I_E
\quad{\rm
and}\\ \o\partial_A\phi &=& 0.
\end{array}
\eeq
See \refsu{ymhymh} for details.  These equations arise by dimensional-reduction
of the 4-dimensional Yang-Mills equations.  Another interpretation is that they
arise if we split projectively flat $SL_2(\C)$-connexions into compact and
non-compact parts (see \refsu{repsta}).

Just as for ordinary Riemann surfaces, the moduli space, $\cal M$, of solutions
to the
\ymh\
equations has an extremely rich geometric structure which we study in the later
sections
of this paper. Let us indicate the main results and outline the contents of
each section.
The first is devoted to preliminaries on orbifold Riemann surfaces and
$V$-bundles (\ie
orbifold bundles):  \refsu{orbint} covers the very basics, for the sake of
revision and in
order to fix notation, while \refsu{orbdiv} deals with the correspondence
between divisors
and holomorphic line $V$-bundles on an orbifold Riemann surface (some of this
may have
been anticipated in unpublished work of B. Calpini).  We particularly draw
attention to
the notational conventions concerning rank 2 $V$-bundles and their rank 1
sub-$V$-bundles
established in \refsu{orbint} which are used throughout this paper.

The second section introduces Higgs $V$-bundles and the appropriate stability
condition (\refsu{highig}) and studies the basic algebraic-geometric properties
of stable Higgs $V$-bundles (\refsu{higalg})---the principal result here is
\refth{stable pairs}.  This material roughly parallels \cite[\S 3]{hi87}, an
important difference being that \cite[proposition 3.4]{hi87} does not
generalise
to the orbifold case.

The third section introduces the \ymh\ equations (\refsu{ymhymh}), discusses
the existence
of solutions on stable Higgs $V$-bundles (\refsu{ymhexi}) and gives the
analytic
construction of $\cal M$ (\refsu{ymhmod}). These first three
subsections parallel \cite[\S\S 4--5]{hi87} and only in \refsu{ymhexi} would
any
significant alteration to Hitchin's work necessary to allow for the orbifold
structure.
The main results are \refth{Narasimhan-Seshadri} and \refth{moduli}.
The Riemannian structure of the moduli space (including the fact that the
moduli space is
hyper-K\"ahler) is also
discussed briefly in \refsu{ymhmod}, following \cite[\S 6]{hi87}. There is one
other
subsection:
\refsu{ymhequ} sketches alternative, equivariant, arguments that can be used
for the
existence theorem and the construction of $\cal M$. This last subsection also
discusses the pull-back map between moduli spaces which arises when an
orbifold Riemann surface is the base of a branched covering by a Riemann
surface---see
\refth{sub}. We stress that equivariant arguments {\em cannot} easily be
applied
throughout the paper---difficulties arise \eg in \refsu{higalg}, \refse{det}
and \refse{rep}.

The fourth section discusses the topology of $\cal M$, following
\cite[\S 7]{hi87}. The results are \refth{Morse} and \refco{topology}.
General formul\ae\  for the Betti numbers are not given
but it is clear how to calculate the Poincar\'e polynomial in any
given instance (however, see \cite{by}).

The fifth section is devoted to the holomorphic symplectic structure on $\cal
M$:  following \cite[\S 8]{hi87}, $\cal M$ is described as a completely
integrable Hamiltonian system via the determinant map $\det :  {\cal M} \to
H^0(K^2)$, defined by taking the determinant of the Higgs field.  This result
is
given as \refth{determinant map} (we believe that a similar result was obtained
by Peter Scheinost).  There are a number of stages to the proof:  first, it is
simpler to use parabolic Higgs bundles and these are discussed in
\refsu{parhig}; \refsu{gendet} contains the major part of the proof, with two
special cases which arise in the orbifold case being dealt with separately in
\refsu{detred} and \refsu{detspe}.  Moreover, it is shown that with respect to
the determinant map $\cal M$ is a fibrewise compactification of the cotangent
bundle of the moduli space of stable $V$-bundles (\refsu{detnon}).

The final section deals with the interpretation of the moduli space as a space
of projectively flat connexions (\refsu{repsta}) or $SL_2(\C)$-representations
of (a central extension of) the orbifold fundamental group (\refsu{reprep}),
the
identification of the submanifold of $SL_2(\R)$-representations
(\refsu{reprea})
and the interpretation of one of the components as Teichm\"uller space
(\refsu{reptei}), which leads to a proof that Teichm\"uller space is
homeomorphic to ball.  The proofs are much like those of \cite[\S\S
9--11]{hi87}
and \cite{do87} and accordingly we concentrate on those aspects of the orbifold
case which are less familiar.

{\em Acknowledgements.} The great debt that the authors owe to the paper
\cite{hi87} is obvious but they are also grateful to Nigel Hitchin for many
useful conversations.  Both authors would also like to thank
Hans Boden, who pointed out an error in \refco{topology}, and Mikio Furuta.

This work is an extension of part of Ben Nasatyr's doctoral thesis:  he would
like to thank Simon Donaldson for his patient supervision and Oscar
Garc\'{\i}a-Prada, Peter Kronheimer and Michael Thaddeus for the contribution
that their comments made to that thesis.  At that time Ben Nasatyr was a
College
Lecturer at Lady Margaret Hall, Oxford, and he spent the following year as a
Post-Doctoral Fellow at the University of British Columbia:  he would like to
thank LMH and NSERC of Canada for their generous financial support and
Gabrielle
Stoy at LMH and David Austin and Dale Rolfsen at UBC for their hospitality.  He
is currently the Sir Michael Sobell Research Fellow at Peterhouse, Cambridge.

Part of Brian Steer's work on this paper took place during a sabbatical year
that he spent in Bonn and Pisa:  he would like to thank Friedrich Hirzebruch
and
the Max-Planck-Institut and Giuseppe Tomassini and the Scuola Normale Superiore
di Pisa for their hospitality.

\bse{Orbifold Riemann Surfaces}{orb}

This section compiles some basic facts about orbifold Riemann surfaces and
fixes some notations
which we will need in the sequel.

\bsu{Introduction to Orbifold Surfaces}{orbint}
We start with the definition and basic properties of orbifold surfaces (or
$V$-surfaces).  The notion of a $V$-manifold was introduced by Satake
\cite{sa56} and re-invented as `orbifold' by Thurston.  By an \de{orbifold
surface} (respectively \de{orbifold Riemann surface}) $M$ we mean a closed,
connected, smooth, real 2-manifold (respectively complex 1-manifold) together
with a finite number (assumed non-zero) of `marked' points with, at each marked
point, an associated order of isotropy $\al$ (an integer greater than one).
(See \cite{sa56} or \cite{sc'83} for full details of the definition.)  Notice
that $M$ has an `underlying' surface where we forget about the marked points
and
orders of isotropy.

Although every point of a surface has a neighbourhood modeled on $D^2$ (the
open
unit disc), we think of a neighbourhood of a {\em marked} point as having the
form $D^2/\Z_\al$, where $\Z_{\al}$ acts on $\R^2 \cong \C$ in the standard way
as the $\al^{\rm th}$ roots of unity.  We make this distinction because $M$ is
to be thought of as an orbifold.  Orbifold ideas do not seem to have been
widely
used in the study of `surfaces with marked points'.  For instance the tangent
$V$-bundle to $D^2/\Z_\al$ is $(D^2 \times \R^2)/\Z_\al$---this leads to an
idea
of an orbifold Riemannian metric on $M$ which corresponds to that of a metric
on
the underlying surface with conical singularities at the marked points (see
\refsu{reptei}).

We introduce the following notations, which will remain fixed throughout this
paper.  Let $M$ be an orbifold (Riemann) surface with topological genus $g$;
denote by $\wo M$ the `underlying' (Riemann) surface obtained by forgetting the
marked points and isotropy.  Denote the number of marked points of $M$ by $n$,
the points themselves by $p_1,\dots,p_n$ and the associated orders of isotropy
by $\al_1,\dots,\al_n$.  Let $\sigma_i$ denote the standard representation of
$\Z_{\al_i}$, with generator $\zeta_i = e^{2\pi\ci/\al_i}$.  At a point where
$M$ is locally $D^2$ or $D^2/\sigma_i$ use $z$ for the standard (holomorphic)
coordinate on $D^2$; call this a local \de{uniformising} coordinate and at a
marked point let $w=z^{\al_i}$ denote the associated local coordinate.  When
giving local arguments centred at a marked point, drop the subscript $i$'s; \ie
use $p$ for $p_i$ and so on.

Given a surface which is the base of a branched covering we naturally consider
it to be an orbifold surface by marking a branch point with isotropy given by
the ramification index.  In this way we arrive at a definition of the
\de{orbifold fundamental group} $\pi_1^V(M)$ (see \cite{sc'83}):  it has the
following presentation
\beql{F-group}
\begin{array}{rcl}\pi_1^V(M) &=& \langle
a_1,b_1,\dots, a_g,b_g,q_1,\dots, q_n \quad |\\ && \quad q_i^{\al_i}=1,\
q_1\dots q_n[a_1,b_1]\dots[a_g,b_g]=1 \rangle.  \end{array}
\eeql
In this presentation $a_1,b_1,\dots,a_g,b_g$ generate the fundamental group of
the underlying surface while $q_1,\dots,q_n$ are represented by small loops
around the marked points.  Similarly, in this situation, the Riemann-Hurwitz
formula suggests the following definition of the \de{Euler characteristic} of
an
orbifold surface:  \beql{Euler characteristic} \chi(M) = 2-2g-n+\sum_{i=1}^n
\frac1{\alpha_i}.\eeql We always work with orbifold surfaces with
$\chi(M)<0$---note that this includes cases with $g=0$ or $g=1$ in contrast to
the situation for ordinary surfaces.

A \de{$V$-bundle}, $E$, with fibre $\C^r$, is as follows.  We ask for a local
trivialisation around each point of $M$ with smooth (or holomorphic) transition
functions; at a marked point $p$ this should be of the form $E|_{D^2/\sigma}
\stackrel{\simeq}{\to} (D^2 \times \C^r)/(\sigma \times \tau)$, where $\tau$ is
an \de{isotropy representation} $\tau :  \Z_\al \to GL_r(\C)$.  We can always
choose coordinates in a $V$-bundle which \de{respect the $V$-structure}:  that
is, if the isotropy representation is $\tau :  \Z_\al \to GL_r(\C)$ then we can
choose coordinates so that $\tau$ decomposes as $\tau = \sigma^{x_1} \oplus
\sigma^{x_2} \oplus \cdots \oplus\sigma^{x_r}$, where, for $j=1,\dots,r$, $x_j$
is an integer with $0\le x_j < \al$ and the $x_j$ are increasing.

We will mostly be interested in rank 2 and rank 1 $V$-bundles and for these we
introduce particular notations for the isotropy, which will be fixed
throughout:
for a rank 2, respectively rank 1, $V$-bundle, denote the isotropy at a marked
point by $x$ and $x'$, respectively by $y$, with $0 \le x,x',y <\al$.  In the
rank 2 case order $x$ and $x'$ so that $x\le x'$.  If a rank 1 $V$-bundle is a
sub-$V$-bundle of a rank 2 $V$-bundle then of course $y \in \{ x,x'\}$:  in
this
case, let $\epsilon \in \{ -1,0,1 \}$ describe the isotropy of the
sub-$V$-bundle, with $\epsilon = 0$ if $x=x'$, $\epsilon = -1$ if $y=x$ and
$\epsilon = 1$ if $y=x'$.  Add subscript $i$'s, when necessary, to indicate the
marked point in question.  Call a vector $(\epsilon_i)$ with $\epsilon_i = 0$
if
$x_i = x'_i$ and $\epsilon_i\in\{\pm1\}$ if not an \de{isotropy vector}.  For a
rank 2 $V$-bundle let $n_0 = \# \{ i :  x_i=x_i'\}$ and for a rank 1
sub-$V$-bundle let $n_\pm = \# \{ i :  \epsilon_i = \pm 1\}$.

If a $V$-bundle is, at a marked point, locally like $(D^2 \times \C^r)/(\sigma
\times \tau)$ then by a Hermitian metric we mean, locally, a Hermitian metric
on
$D^2 \times \C^r$ which is equivariant with respect to the action of $\Z_\al$
via $\sigma\times\tau$.  Considering the tangent $V$-bundle, we can also define
the concepts of Riemannian metric and orientation for an orbifold surface (an
orientation of an orbifold surface is just an orientation of the underlying
surface).

We introduce the notion of a connexion in a $V$-bundle in the obvious way.  The
first Chern class or degree of a $V$-bundle can be defined using Chern-Weil
theory.  Notice that the degree of a $V$-bundle is a {\em rational} number,
congruent modulo the integers to the sum $\sum_{i=1}^n(y_i/\al_i)$, where
$(y_i)$ is the isotropy of the determinant line $V$-bundle.

When $E$ is a rank 2 $V$-bundle with isotropy $(x_i,x_i')$, as above then we
write
\beq
c_1(\Lambda^2E) = l + \isum\frac{x'_i + x_i}{\al_i},
\eeq
for $l \in \Z$. Similarly, if $L$ is a sub-$V$-bundle with isotropy given by an
isotropy
vector $(\epsilon_i)$ in the manner explained above then we write
\beq
c_1(L) = m + \isum \frac{\epsilon_i(x'_i - x_i) + (x'_i + x_i)}{2\al_i}
\eeq
for $m \in \Z$. These meanings of $l$ and $m$ will be fixed throughout.

Topologically, $U(1)$ and $U(2)$ $V$-bundles are classified by their isotropy
representations and first Chern class:  we quote the following classification
result from \cite{fs92}.
\bprn{Furuta-Steer}{v-bundles} Let $M$ be an orbifold surface.  Then, over $M$:
\begin{enumerate}
\item any complex line $V$-bundle is topologically determined by its isotropy
representations and
degree, \item any $SU(2)$ $V$-bundle is topologically determined by its
isotropy representations
(necessarily of the form $\sigma^{x}\oplus\sigma^{-x}$, where $0\le x\le
[\al/2]$) and \item any
$U(2)$ $V$-bundle is topologically determined by its isotropy representations
and its determinant
line $V$-bundle.  \end{enumerate} \epr
\bre{subbundles}
Let $E$ be a $U(2)$ $V$-bundle with isotropy $(x_i,x_i')$ and let
$(\epsilon_i)$ be any
isotropy vector. Then there exists
a $U(1)$ $V$-bundle $L$ with isotropy specified by $(\epsilon_i)$ (unique up to
twisting
by a $U(1)$-bundle \ie up to specifying the integer $m$, above) and,
topologically,
$E=L\oplus L^*\Lambda^2E$, by
\refpr{v-bundles}.
\ere

\bsu{Divisors and Line $V$-bundles}{orbdiv}

The theory of divisors developed here has also been dealt with in the Geneva
dissertation of B. Calpini written some time ago.

Suppose $M$ is an orbifold Riemann surface.  It is convenient to associate an
order of isotropy $\al_p$ to every point $p$; it is 1 if the point is not one
of
the marked points (and $\al_i$ if $p=p_i$ for some $i$).  A \de{divisor} is
then
a linear combination \beq D = \sum_{p\in M}\frac{n_p}{\al_p}.p \eeq with
$n_p\in
\Z$ and zero for all but a finite number of $p$.

If $f$ is a non-zero meromorphic function on $M$ we define the \de{divisor of
$f$} by $Df = \sum_p \nu_p(f).p$.  Here $\nu_p(f)$ is defined in the usual way
when $\al_p=1$.  When $\al_p=\al > 1$ and $z$ is a local uniformising
coordinate
with $\rho :  D^2 \surjarrow D^2/\sigma$ the associated projection, then on
$D^2$ we find that $\rho^*f$ has a Laurent expansion of the form \beq
\sum_{j\ge
-N}a_j z^{\al j}\qquad\mbox{with $a_{-N}\ne 0$} \eeq and we set $\nu_p(f) =
-N$.
(The divisor of a meromorphic function is thus an {\em integral} divisor.)  Two
divisors $D$ and $D'$ are \de{linearly equivalent} if \beq D-D' = Df \eeq for
some meromorphic f. The \de{degree} of a divisor $D=\sum_p(n_p/\al_p).p$ is
defined to be $d(D)=\sum_p n_p/\al_p$.

The correspondence between divisors and holomorphic line $V$-bundles goes
through in exactly the same way as for Riemann surfaces without marked points.
To a point $p$ with $\al_p=1$ we associate the point line bundle $L_p$ as in
\cite{gu66}.  If $\al_p=\al>1$ then to the divisor $p/\al$ we associate the
following $V$-bundle.  Let $z$ be a local uniformising coordinate; then, making
the appropriate identification locally with $D^2/\sigma$, we define \beq
L_{p/\al} = ((D^2\times \C)/(\sigma\times\sigma)) \cup_\Phi
((M\setminus\{p\})\times\C ), \eeq where $\Phi:  (D^2\setminus\{0\}\times
\C)/(\sigma\times\sigma) \to ((M\setminus\{p\})\times\C )$ is given by its
$\Z_\al$-equivariant lifting \beq \lift\Phi:(D^2\setminus\{0\})\times\C &\to&
((D^2/\sigma)\setminus\{0\})\times\C\\ (z,z')&\mapsto&(z^{\al},z^{-1}z').  \eeq

This $V$-bundle has an obvious section `$z$'; this is given on $D^2\times \C$
by
$z\mapsto(z,z)$ and extends by the constant map to the whole of $M$.  So
$L_{p/\al}$ is positive.  We denote by $L_i$ the line $V$-bundle
$L_{p_i/\al_i}$, associated to the divisor $p_i/\al_i$, and by $s_i$ the
canonical section `$z$'.

Finally for a general divisor \beq D = \sum_{p\in M}\frac{n_p}{\al_p}.p \eeq we
set \beq L_D = \otimes_p(L_{p/\al_p})^{n_p}.  \eeq

As for a meromorphic function, we can define the divisor of a meromorphic
section of a line $V$-bundle $L$.  If $p$ has ramification index $\al_p=\al$
and
we have a local uniformising coordinate $z$ and a corresponding local
trivialisation $L|_{D^2/\sigma} \cong (D^2\times\C) /(\sigma\times\sigma^y) $,
for some isotropy $y$ (with, by convention, $0\le y < \al$), then locally we
have $ s(z) = \sum_{j\ge -N'}a'_j z^j$ with $a'_{-N'}\ne 0$.  However, we have
$\Z_\al$-equivariance which means that $s(\zeta.z)=\zeta^ys(z)$ (where $\zeta =
e^{2\pi\ci/\al}$ generates $\Z_\al$).  It follows that $a'_j = 0$ unless
$j\equiv y\pmod \al$ and hence \beql{Taylor} s(z) = z^{y}\sum_{j\ge -N}a_j
z^{\al j}\qquad\mbox{with $a_{-N}\ne 0$,} \eeql where $-N\al + y = -N'$.  We
define $\nu_p(s)=-N'/\al = -N + y/\al$:  so for the canonical section $s_i$ of
the line $V$-bundle $L_{i}$ we have $\nu_{p_i}(s_i) = 1/\al_i$.

\bpr{divisors} The above describes a bijective correspondence between
equivalence classes of
divisors and of holomorphic line $V$-bundles.  The degree $d(D)$ of a divisor
$D$ is just $c_1(L_D)$, the first Chern class of the corresponding line
$V$-bundle.  \epr
\bpf Much of the proof is contained in \cite{fs92}.  The correspondence has
been
defined above and it is clear that if we start from a divisor $D$ and pass to
$L_D$ then taking the divisor associated to the tensor product of the canonical
sections we get back $D$.  We have to show that the correspondence behaves well
with respect to equivalence classes.  If $D_1\equiv D_2$, where $D_j = \sum
(n^{(j)}_p/\al_p).p$ for $j=1,2$, then from what we know about divisors of
meromorphic functions we see that $n^{(1)}_p \equiv n_p^{(2)} \pmod{\al_p}$.
Now $L_{D_j}=\otimes_p(L_{p/\al_p})^{n^{(j)}_p}$.  Since $n^{(1)}_p \equiv
n_p^{(2)} \pmod{\al_p}$, we find that
$L_{D_j}\otimes\bigotimes_{i=1}^n(L_{p_i/\al_i})^{-n^{(1)}_{p_i}}$ is a genuine
line bundle for $j=1,2$.  Moreover the two are equivalent because the
corresponding divisors are.  Hence $L_{D_1}\equiv L_{D_2}$.  Similarly we show
that two meromorphic sections of the same line \vb\ define equivalent divisors.
\epf

\bco{divisors1} If
$L$ is a
holomorphic line $V$-bundle with $c_1(L)\le 0$ then $H^0(L) = 0$, unless $L$ is
trivial.  \eco

Let $L$ be a holomorphic line $V$-bundle over $M$, with isotropy $y_i$ at
$p_i$,
and let $\co(L)$ be the associated sheaf of germs of holomorphic sections; we
take the cohomology of $L$ over $M$ to be the sheaf cohomology of $\co(L)$ over
$\wo M$.  From \refeq{Taylor}, $\co(L)$ is locally free over $\co_M = \co_{\wo
M}$ and hence there is a natural line bundle $\wo L$ over $\wo M$ with $\co(\wo
L) \cong \co(L)$.  If we define $\wo L = L\otimes L_1^{-y_1}\otimes
\cdots\otimes L_n^{-y_n}$ then this gives the required isomorphism of sheaves.

\bpr{parabolic} If $L$ is a holomorphic line $V$-bundle then, with $\wo L$
defined as above, there
is a natural isomorphism of sheaves $\co(L)\cong \co(\wo L)$ given by tensoring
with the canonical
sections of the $L_i$.
\epr
\bpf
Recall that $\ilist{s}$ are the canonical sections of $\ilist{L}$.  If $s$ is a
holomorphic section of $L$ then $\wo s=s_1^{-y_1}\dots s_n^{-y_n}s$ will be a
meromorphic section of $\wo L$, holomorphic save perhaps at $p_i$.  In fact (by
choosing a local coordinate) we see that $\wo s$ has removable singularities at
$p_i$ and that $D(\wo s) = Ds - \isum ({y_i}/{\al_i})p_i$.  Conversely, given a
section $\wo s$ of $\wo L$, then $s_1^{y_1}\dots s_n^{y_n}\wo s$ is a section
of
$L$ and the correspondence is bijective.
\epf

As corollaries we get the orbifold Riemann-Roch theorem, originally due to
Kawasaki \cite{ka79} and an orbifold version of Serre duality.
\bthn{Kawasaki-Riemann-Roch}{Riemann-Roch} Let $L$ be a
holomorphic line $V$-bundle
with the isotropy at $p_i$ given by ${y_i}$, with $0 \le y_i < \al_i$.  Then $$
h^0(L) - h^1(L)
= 1-g+c_1(L) -
\sum_{i=1}^n \frac {y_i}{\alpha_i}, $$ where $h^i$ denotes the dimension of
$H^i$.  \eth

\bth{Serre duality} If $L$ is a holomorphic line \vb\ and $K_M$ is the
canonical
$V$-bundle of the
orbifold Riemann surface then \beq H^1(L) \cong H^{0}(L^* K_M)^*.  \eeq \eth
\bpf By
definition, $H^1(L)=H^1(\co(L))=H^1(\wo L)$.  So $H^1(L)\cong H^{0}((\wo L)^*
K_{\wo M})^*$ by the standard duality.  But $(\wo L)^* K_{\wo M}=\wo {L^* K_M}$
by
a straightforward computation.  \epf

\bse{Higgs $V$-Bundles}{hig}

Throughout this section $E \to M$ is a holomorphic rank 2 $V$-bundle
over an orbifold Riemann surface with $\chi(M)<0$ and we write $K=K_M$, the
canonical $V$-bundle,
and $\Lambda = \Lambda^2E$, the determinant line $V$-bundle.

\bsu{Higgs $V$-Bundles}{highig}

In this subsection we introduce Higgs $V$-bundles---this is a straightforward
extension of the basic
material in Hitchin's paper \cite{hi87} to orbifold Riemann surfaces.

Define a \de{Higgs field}, $\phi$, to be a
holomorphic section of ${\rm End}_0(E)\otimes K$ where
${\rm End}_0(E)$ denotes the trace-free endomorphisms of $E$.  A \de{Higgs
$V$-bundle} or \de{Higgs
pair} is just a pair $({E},\phi)$.

Let $({E}_1,\phi_1)$ and $({E}_2,\phi_2)$ be two Higgs $V$-bundles.  A
\de{homomorphism of Higgs
$V$-bundles} is just a homomorphism of $V$-bundles $h :  E_1 \to E_2$ such that
$h$ is holomorphic
and intertwines $\phi_1$ and $\phi_2$.  The corresponding notion of an
\de{isomorphism of Higgs
$V$-bundles} is then clear.

A holomorphic line sub-$V$-bundle $L$ of $E$ is called a \de{Higgs
sub-$V$-bundle} (or
`$\phi$-invariant sub-$V$-bundle') if $\phi(L) \subseteq KL$.  A Higgs
$V$-bundle $({E},\phi)$ is
said to be \de{stable} if \beql{stable Higgs} c_1(L) < \frac12
c_1(E),\quad\mbox{for every
Higgs sub-$V$-bundle, $L$.} \eeql If we allow possible equality in
\refeq{stable Higgs} then the
Higgs $V$-bundle is called \de{semi-stable}.  If a Higgs $V$-bundle is stable
or a direct sum of two
line $V$-bundles of equal degree with $\phi$ also decomposable then (it is
certainly
semi-stable and) it is called
\de{polystable}. If $E$ is stable then certainly $(E,\phi)$ is stable for
any Higgs field $\phi$. The following result, due to Hitchin in the smooth
case
\cite[proposition 3.15]{hi87}, goes over
immediately.
\bpr{stable regular} Let $({E}_1,\phi_1)$ and $({E}_2,\phi_2)$ be stable Higgs
$V$-bundles with isomorphic holomorphic determinant line $V$-bundles,
$\Lambda^2{E}_1 \cong
\Lambda^2{E}_2$.  Suppose that $\psi:  E_1 \to E_2$ is a non-zero homomorphism
of Higgs $V$-bundles.
Then $\psi$ is an isomorphism of Higgs $V$-bundles.  If
$({E}_1,\phi_1)=({E}_2,\phi_2)$ then $\psi$
is scalar multiplication.  \epr

\bsu{Algebraic Geometry of Stable Higgs $V$-Bundles}{higalg}

For applications in later sections we now develop some results on the
possibilities for
stable Higgs $V$-bundles. Higgs $V$-bundles are
holomorphic $V$-bundles with an associated `Higgs field'; a holomorphic
$(1,0)$-form-valued
endomorphism of the $V$-bundle.  We assume familiarity with \cite[\S
3]{hi87}.

Given $E \to M$, we investigate whether there are any Higgs fields $\phi$ such
that the Higgs pair
$(E,\phi)$ is stable.  Recall that the isotropy of $E$ at $p_i$ is denoted by
$(x_i,x_i')$ and that
$n_0=\#\{i\,:\,x_i = x_i'\}$.  We will suppose throughout that $n_0 < n$---this
is because the case
$n=n_0$ is just that of a genuine bundle twisted by a line $V$-bundle and so
essentially
uninteresting (see also \refsu{detred}).

The following lemma is a simple computation using the
Kawasaki-Riemann-Roch theorem and Serre duality.
\ble{}
We have \[ h^0(K^2) = \chi(K^2) = 3g-3+n \and \chi(\End_0(E)\otimes
K)=3g-3+n-n_0. \] \ele

If $E$ is stable we know that the only endomorphisms of $E$ are scalars and so
$h^0(\End_0(E))=0$;
consequently if $3-3g-n+n_0>0$ (this only happens if $g=0$ and $n-n_0\le 2$)
there are no stable
$V$-bundles.

Suppose that $L$ is a holomorphic sub-$V$-bundle of $E$.  Then we have the
short exact sequences \beql{b}
&0\to L \stackrel{i}{\to} E \stackrel{j}{\to} L^* \Lambda\to 0&\and\nonumber\\
&0\to
L  \Lambda^* \stackrel{j^*}{\to} E^* \stackrel{i^*}{\to} L^*\to 0& \eeql from
which follows
\beql{d} 0\to E^*\otimes KL \to \End_0(E) \otimes K\to KL^{-2}  \Lambda\to 0.
\eeql
Associated to \refeq{b} tensored by $KL$ is the long exact sequence in
cohomology \beql{lesb}
\begin{array}{c} 0 \to H^0(KL^2\Lambda^*) \to H^0(E^*\otimes KL) \to
H^0(K)\stackrel{\delta}{\to}\qquad\\ \qquad\stackrel{\delta}{\to}
H^1(KL^2\Lambda^*) \to
H^1(E^*\otimes KL) \to H^1(K) \to 0 \end{array} \eeql and associated to
\refeq{d} we have \beql{les}
\begin{array}{c} 0 \to H^0(E^*\otimes KL) \to H^0(\End_0(E) \otimes K) \to
H^0(KL^{-2}
\Lambda)\stackrel{\delta}{\to}\qquad\\ \qquad\stackrel{\delta}{\to}
H^1(E^*\otimes KL) \to
H^1(\End_0(E) \otimes K) \to H^1(KL^{-2}  \Lambda) \to 0. \end{array} \eeql

Now let us review the strategy of the proof of \cite[proposition 3.3]{hi87}:
if $E$ is stable then
all pairs $(E,\phi)$ are certainly stable and we know something about stable
$V$-bundles from
\cite{fs92}.  If $E$ is not stable then there is a destabilising sub-$V$-bundle
$L_E$.  Recall that $L_E$
is unique if $E$ is not semi-stable.  Moreover, in the semi-stable case the
assumption $n\ne n_0$
implies that $L_E \not\cong L_E^*\Lambda$ and so $L_E$ is unique if $E$ is not
decomposable and
if it is then $L_E$ and $L_E^*\Lambda$ are the only destabilising
sub-$V$-bundles.  Thus there will be
some $\phi$ such that the pair $(E,\phi)$ is stable unless every Higgs field
fixes $L_E$ (or
$L_E^*\Lambda$, in the semi-stable, decomposable case).
Moreover, the subspace of sections leaving $L$ invariant is $H^0(E^*\otimes
KL)
\subset
H^0(\End_0(E)\otimes K)$. It follows that a necessary and sufficient condition
for $E$ to
occur in a stable pair is $H^0(E^*\otimes KL_E)\ne H^0(\End_0(E)\otimes K)$
(and similarly for
$L_E^*\Lambda$, in the semi-stable, decomposable case).  Considering
\refeq{les} this amounts to
non-injectivity of the Bockstein operator $\delta$, which we consider in the
next
lemma---proved as in the proof of \cite[proposition 3.3]{hi87}. From the lemma
we obtain a version of \cite[proposition 3.3]{hi87}.
\ble{extension'} If $L$ is a
sub-$V$-bundle of $E$ with $\deg(L)\ge \deg(\Lambda)/2$ then \begin{enumerate}
\item\label{ex1}
$H^1(E^*\otimes KL)\cong\C;$ \item\label{ex2} $H^0(KL^{-2}
\Lambda)\stackrel{\delta}{\to}H^1(E^*\otimes KL)$ is surjective if and only if
$e_E \ne 0$,
where $e_E\in H^1(L^2  \Lambda^*)$ is the extension class.  \end{enumerate}
\ele \bpf
\begin{enumerate} \item Consider the long exact sequence in cohomology
\refeq{lesb} for $L$, which
includes the segment \beql{bit} \cdots \to H^1(KL^2  \Lambda^*)
\stackrel{j^*}{\to}
H^1(E^*\otimes KL) \stackrel{i^*}{\to} \C \to 0 . \eeql Then the result follows
from the fact
that $h^1(KL^2  \Lambda^*)=0$, using Serre duality and the vanishing theorem.
\item
Consider \refeq{les} and let $i^*$ be the map on cohomology indicated in
\refeq{bit}; then the
result follows from the fact that $i^*.\delta$ is multiplication by the
extension class $e_E$.
\end{enumerate} \epf

\bpr{non-s} Let $E$ be a non-stable $V$-bundle.  Then $E$ appears
in a stable pair if and only if one of the following holds:  \begin{enumerate}
\item\label{ns1} $E$
is indecomposable with $h^0(KL_E^{-2}\Lambda)>1$; \item\label{ns2} $E$ is
decomposable, not
semi-stable with $h^0(KL_E^{-2}\Lambda)\ge 1$; \item\label{ns3} $E$ is
decomposable,
semi-stable with $h^0(KL_E^{-2}\Lambda)\ge 1$ and $h^0(KL_E^{2}\Lambda^*)\ge
1$.
\end{enumerate} \epr

To find more precise results in the case that $E$ is semi-stable we
estimate $h^0(KL_E^{-2}\Lambda)$ and $h^0(KL_E^{2}\Lambda^*)$ using the
following lemmas. For these recall the definitions of the integers
$n_0$, $n_\pm$, $l$ and $m$ from \refsu{orbint}.
\ble{chi}
Suppose that $L$ is any sub-$V$-bundle of $E$.  Then, with the notations
established
above,
\beq \chi(KL^{-2}  \Lambda) &=& l-2m
+g - 1 + n_-\and\\
\chi(KL^{2}  \Lambda^*)
&=& 2m-l +g - 1  +  n_+. \eeq  Moreover:
\begin{enumerate}
\item\label{xpos} if $2c_1(L)-c_1(\Lambda) \ge 0$ then $h^0(KL^{2}  \Lambda^*)
=
\chi(KL^{2}  \Lambda^*) \ge g$ and $\chi(KL^{-2}  \Lambda) \le g - 2 + n -
n_0$;
\item\label{xneg} if $2c_1(L)-c_1(\Lambda) \le 0$ then
$h^0(KL^{-2}  \Lambda) = \chi(KL^{-2}  \Lambda) \ge g$ and
$\chi(KL^{2}  \Lambda^*) \le g - 2 + n - n_0$.
\end{enumerate}
\ele
\bpf The first part is just the Kawasaki-Riemann-Roch theorem. Now consider
\refpa{xpos}
(\refpa{xneg} is entirely similar):
we have $H^1(KL^{2}\Lambda^*) \cong H^0(L^{-2}\Lambda)^*$ and this is zero
(because the degree is
non-positive and the isotropy is non-trivial as $n> n_0$).
Let $\theta = \sum_{i=1}^{n}{\epsilon_i(x'_i-x_i)}/{\al_i}$ so that
$2c_1(L)-c_1(\Lambda)
\equiv \theta \pmod{\Z}$. Then $-n_- < \theta < n_+$ and the estimates on
$\chi(KL^{2}
\Lambda^*)$ and $\chi(KL^{-2}
\Lambda)$ follow.
\epf

\ble{bounds}
For a given $M$ and $n-n_0$, an $E$ (with the given $n-n_0$) such that the
bounds on
$\chi(KL^{2}
\Lambda^*)$ and $\chi(KL^{-2}  \Lambda)$ in \refle{chi}, parts 1 and 2
are attained exists if and only if
\beq
\min_{ \{i_1,\dots,i_{n-n_0}\} \subseteq \{1,\dots,n \} }\left\{
\sum_{j=1}^{n-n_0} \frac{1}{\al_{i_j}} \right\} \le 1.\eeq
For a given topological $E$ the bounds are attained for some
holomorphic structure on $E$ if and only if
\beq
\min_{ \{(\epsilon_i)\  :\  n_+ + l\equiv 1 (2) \}} \left\{ n_+  -
\isum\frac{\epsilon_i(x'_i-x_i)}{\al_i} \right\} \le 1,\eeq
where $(\epsilon_i)$ varies over all isotropy vectors with $n_+ +l \equiv 1
(2)$.
\ele
\bpf
To see this we construct examples as follows.  It is sufficient to consider
only
topological examples and therefore, given any $M$ and topological $E$, to
choose
$(\epsilon_i)$ and $m\in \Z$ to specify $L$ topologically.  (Examples where $L$
is a topological sub-$V$-bundle of $E$ exist by \refre{subbundles}.)

Now, given a choice of $(\epsilon_i)$ and $m$, we
have $\chi(KL^{-2}\La) = l - 2m + g - 1 + n_-$ and $\chi(KL^{2}\La^*) = 2m - l
+ g - 1 +
n_+$ from \refle{chi}.
So, for $2c_1(L)-c_1(\Lambda) \ge 0$ (the case $2c_1(L)-c_1(\Lambda)\le 0$ is
entirely similar), the bounds
are attained provided $2m - l + n_+ = 1$ and $2m - l +
\sum_{i=1}^{n}\epsilon_i(x'_i-x_i)/\al_i \ge 0$.
Since we can vary $m$, the first equation just fixes the parity of $n_+$.
Hence the problem reduces to finding $(\epsilon_i)$ such that
\beql{succinct}
\sum_{i=1}^{n}\frac{\epsilon_i(x'_i-x_i)}{\al_i} - n_+ &\ge& -1\and\\
n_+ + l &\equiv& 1 \pmod{2}.\label{eq:constraint}
\eeql
This gives the desired result, for a given topological $E$. To see whether
examples exist for a given $M$ and $n-n_0$ as we allow $E$ to vary over
topological types with fixed $n-n_0$, we simply note that the maximum value of
the left-hand side of \refeq{succinct} (subject to \refeq{constraint}) is \beq
\max_{ \{i_1,\dots,i_{n-n_0}\} \subseteq \{1,\dots,n \} }\left\{
\sum_{j=1}^{n-n_0} \left( -\frac{1}{\al_{i_j}} \right) \right\}.\eeq
Thus the bounds are certainly attained if the $\al_i$ are such that this is not
less
than $-1$.\epf

\bco{semi-stable'} If $L$ is a
sub-$V$-bundle of $E$ with
$c_1(L)=c_1(\Lambda)/2$ and $\epsilon_i$, $n_+$ and $n_-$ are defined by the
isotropy of $L$, as
before, then \beql{} h^0(\End_0(E)\otimes K)
&=& \left\{\begin{array}{ll} 3g-3+n-n_0 & \mbox{ if }0\to L \to E \to L^*
\Lambda\to
0\mbox{ is non-trivial;}\\ 3g-2+n-n_0 & \mbox{ if it is trivial;}
\end{array}\right.\label{h0end}\\
h^0(E^*\otimes KL) &=& 2g - 1 -
\sum_{i=1}^{n}\frac{\epsilon_i(x'_i-x_i)}{\al_i} + n_+;
\label{vkl}\\
h^0(KL^{-2}  \Lambda) &=& g - 1 + \sum_{i=1}^{n}
\frac{\epsilon_i(x'_i-x_i)}{\al_i} +
n_-;\label{kl}\\
h^0(KL^{2}  \Lambda^*) &=& g - 1 - \sum_{i=1}^{n}
\frac{\epsilon_i(x'_i-x_i)}{\al_i}
+ n_+;\nonumber.  \eeql
Moreover,
\beq 2g \le &h^0(E^*\otimes KL)& \le
n-n_0 + 2g -2,\\ g \le &h^0(KL^{-2}  \Lambda) & \le n-n_0 + g -2\and\\ g \le
&h^0(KL^{2}  \Lambda^*) & \le n-n_0 + g -2.  \eeq
These estimates are attained for all values
of $g$ and $n-n_0$ (but not necessarily for all $M$ or $E$).
\eco \bpf The results on $h^0(KL^{-2}  \Lambda)$ and
$h^0(KL^{2}  \Lambda^*)$  follow from \refle{chi}.  Moreover we know that
$h^1(E^*\otimes KL)= 1$
from \refle{extension'}, \refpa{ex1} and so $h^0(E^*\otimes KL)$ follows from
the
Kawasaki-Riemann-Roch theorem. To calculate $h^0(\End_0(E)\otimes K)$ we use
\refeq{les} and \refle{extension'}, \refpa{ex2}.  The estimates on
$h^0(KL^{-2}\Lambda)$ and $h^0(KL^2\Lambda^*)$ are contained
in \refle{chi} and the estimate on $h^0(E^*\otimes KL)$ follows (as
$h^0(E^*\otimes KL) = - h^0(KL^{-2}  \Lambda) +
3g -2 + n - n_0$). \epf
When $c_1(L)=c_1(\Lambda)/2$ it is not possible to have $n-n_0=1$
(because $c_1(L^2\Lambda^*)$ cannot be an integer if $n-n_0=1$ but, on the
other hand, it
is supposed zero).

Applying these results to $L_E$ (and $L_E^*\Lambda$ in the semi-stable,
decomposable case)
we can strengthen \refpr{non-s} as far as it refers to semi-stable
$V$-bundles. Adding in some necessary conditions on $g$ and $n-n_0$ derived
from
our estimates above we obtain the following theorem.
\bth{stable pairs} A holomorphic rank 2
$V$-bundle $E$ occurs in a stable pair if and only if one of the following
holds:  \begin{enumerate}
\item\label{s} $E$ is stable (if $g=0$ then necessarily $n-n_0\ge 3$);
\item\label{ss} $E$ is semi-stable, not
stable (necessarily $n-n_0 \ge 2$) with one of the following holding:
\begin{enumerate}
\item\label{ssin1} $E$ is indecomposable and $g>1$;
\item\label{ssin2} $E$ is indecomposable, $g=0$
or 1 and $h^0(KL_E^{-2}  \Lambda)>1$ (necessarily $g + n-n_0\ge 4$);
\item\label{ssde1} $E$
is decomposable and $g>0$;
\item\label{ssde2} $E$ is decomposable, $g=0$ and $1 \le
h^0(KL_E^{-2}  \Lambda) \le n -n_0 -3$ (necessarily $n-n_0\ge 4$);
\end{enumerate}
\item\label{not semi-stable}
$E$ is not semi-stable with one of the following holding:  \begin{enumerate}
\item\label{nsin} $E$
is indecomposable and $h^0(KL_E^{-2}  \Lambda)>1$ (necessarily $g\ge 2$ or $g +
n-n_0 \ge
4$; if $g=2$ and $n-n_0 =1$ then $\wo{KL_E^{-2} \Lambda}$ is necessarily
canonical);
\item\label{nsde} $E$ is decomposable and $h^0(KL_E^{-2}
\Lambda) \ge 1$
(necessarily $g\ge 1$ or $n-n_0\ge 3$; if $2g+n-n_0=3$ then $\wo{KL_E^{-2}
\Lambda}$ is necessarily
trivial).  \end{enumerate} \end{enumerate}
In all cases the necessary conditions are the best possible ones depending only
on $g$ and
$n-n_0$. \eth
\bpf
In \refpa{ss} the first three items follow from
\refco{semi-stable'} together with \refpr{non-s}, parts 1 and 3,
while for the last item we note that when $g=0$, $h^0(KL_E^{2}  \Lambda^*)\ge
1$
if and only if $h^0(KL_E^{-2}  \Lambda)<n-n_0-2$ (from \refco{semi-stable'})
and
apply \refpr{non-s}, \refpa{ns3}.

Only the necessary conditions in \refpa{not semi-stable} need any additional
comment.  Using \refle{chi}, \refpa{xpos} we have that $\chi(KL_E^{-2} \Lambda)
\le g - 2 + n - n_0$ and the bound is attained for some $M$ and $E$ by
\refle{bounds}.  Thus if $g>2$ there are cases with $\chi(KL_E^{-2} \Lambda)\ge
2$ and hence $h^0(KL_E^{-2} \Lambda)\ge 2$.  If $g=2$ then there are cases with
$\chi(KL_E^{-2} \Lambda)=n-n_0$, similarly.  The only problem then occurs if
$n-n_0 =1$ when $c_1(\wo{KL_E^{-2} \Lambda})=2$:  in order to have
$h^0(KL_E^{-2} \Lambda)>1$ we must have $\wo{KL_E^{-2} \Lambda} = K_{\wo M}$.

Similarly, if $g=1$ we can suppose that $\chi(KL_E^{-2} \Lambda) = n - n_0 -1$.
Then for $h^0(KL_E^{-2}\Lambda) >1$ we need $n-n_0 \ge 3$ and for
$h^0(KL_E^{-2}\Lambda) \ge 1$ we need $n-n_0 \ge 1$ with $\wo{KL_E^{-2}
\Lambda}$ trivial if $n-n_0 =1$.  Finally, if $g=0$ we need $n-n_0 \ge 4$ for
$h^0(KL_E^{-2}\Lambda) >1$ and $n-n_0 \ge 3$ (with $\wo{KL_E^{-2} \Lambda}$
trivial if $n-n_0 =3$) for $h^0(KL_E^{-2}\Lambda) \ge 1$.
\epf
For each of the items of \refth{stable pairs} examples of such
$V$-bundles do
actually exist (see also \refse{top} and \refse{det}).  Only items \ref{ssin2},
\ref{ssde2},
\ref{nsin} and \ref{nsde} pose any problem but it
is fairly easy to construct the required examples using the ideas of
\refsu{orbdiv} and
\refle{bounds}. Of particular
interest is \refpa{nsde} when $g=0$ and $n-n_0=3$:  we have the following
result (compare
\refse{top}).
\bpr{g0} There exist orbifold Riemann surfaces with $g=0$ with $V$-bundles
with
$n-n_0=3$ over them which are decomposable but not semi-stable and
exist in stable pairs. Such a stable pair contributes an isolated point
to the moduli space (which is nevertheless connected---see \refco{topology}).
\epr
\bpf
We set $E=L_E\oplus L_E^*\Lambda$ with $2c_1(L_E) >
c_1(\Lambda)$. Now, according
to \refth{stable pairs}, \refpa{nsde}, we get a stable pair if and only if
$\wo{KL_E^{-2}  \Lambda}$ is trivial. Moreover, applying \refsu{orbdiv} or
\refle{bounds}, we see
that examples certainly exist.

We write the Higgs field according to the decomposition $\phi = \left(
\begin{array}{cc}t & u\\ v &
-t \end{array}\right)$.  Now $h^0(KL_E^{-2} \Lambda)=1$ implies that
$h^0(KL_E^2  \Lambda^*)=0$ and hence $u=0$. More simply, $g=0$ implies $t=0$
and so $\phi$ is given by $v$, with
$v\in H^0(KL_E^{-2}  \Lambda)\cong \C$ non-zero for a stable pair.  Now we
need
to consider the
action of $V$-bundle automorphisms: $\left( \begin{array}{cc}\lambda & 0\\ 0 &
\lambda^{-1} \end{array}\right)$ acts on $ H^0(KL_E^{-2}  \Lambda)\cong \C$ by
$z \mapsto
\lambda^2 z$ and hence there is a single orbit.
\epf

Notice that \cite[proposition 3.4]{hi87} does {\em not} extend to orbifold
Riemann surfaces with $\chi(M)<0$.  To prove that result Hitchin uses
Bertini's
theorem to show that, for
a given rank 2 holomorphic bundle over a Riemann surface with negative Euler
characteristic, either
the generic Higgs field leaves no subbundle invariant or there is a subbundle
invariant under all
Higgs fields; he then shows that the latter cannot happen when the bundle
exists in a stable pair.
Although we have not been able to enumerate all the cases in which this result
is false in the
orbifold case there are three things which can go wrong:
\begin{enumerate}
\item Bertini's theorem may not apply and the conclusion may be false: $E$ may
be such that it
exists in a stable pair, the generic Higgs field has an invariant
sub-$V$-bundle and no
sub-$V$-bundle is invariant by all Higgs fields;
\item $E$ may be stable and have a sub-$V$-bundle invariant by
all Higgs fields;
\item $E$ may be non-stable, exist in a stable pair and have a
sub-$V$-bundle invariant by all Higgs fields.
\end{enumerate}
We give counterexamples of the first and third types.  Although we suspect that
counterexamples of the second type also exist we have not been able to show
this.  For a counterexample where Bertini's theorem doesn't apply consider the
following:  if $g=1$ and $n-n_0=1$ then, anticipating \refle{invariants}, {\em
every Higgs field has an invariant sub-$V$-bundle} and yet if $E$ is a
non-stable $V$-bundle which exists in a stable pair (these exist by
\refth{stable pairs}, \refpa{nsde}) then {\em no sub-$V$-bundle is invariant by
all Higgs fields.} All counterexamples of the third type are given in the
following proposition, which also has interesting applications in \refse{top}.

\bpr{counter}
A non-stable $V$-bundle $E$ exists in a stable pair and has a sub-$V$-bundle
invariant by all Higgs
fields if and only if $g=0$, $E=L_E \oplus L_E^*\La$ with $2c_1(L_E) >
c_1(\La)$ and $L_E$ is
such that the bounds in \refle{chi}, \refpa{xpos} are attained. Moreover, there
exist orbifold
Riemann surfaces with such $E$ over them, with $E$ having any given $n-n_0 \ge
3$.
\epr
\bpf
Suppose $E$ is non-stable, exists in a stable pair and has a sub-$V$-bundle
invariant by all Higgs
fields. Since $E$ is non-stable and exists in a stable pair the destabilising
sub-$V$-bundle(s)
cannot be invariant by all Higgs fields. Moreover, if $h^0(KL_E^2\La^*)>0$
then, via the
inclusions $H^0(KL_E^2\La^*) \hookrightarrow H^0(E^*\otimes KL_E)
\hookrightarrow
H^0(\End_0(E)\otimes K)$, we get a family of Higgs fields which leave no
sub-$V$-bundle except $L_E$
invariant---hence we must have $h^0(KL_E^2\La^*)=0$. By \refle{chi},
\refpa{xpos} this can only
happen if the bounds there are attained and $g=0$. Now consideration of the
long exact sequence
\refeq{lesb}
shows that $h^0(E^*\otimes KL_E)=g$ and hence \refle{extension'} and
\refeq{les} together show that
$E$ is decomposable. Considering the Higgs field according to the
decomposition, in the manner of
\refpr{g0}, we see that $L_E^*\La$ is invariant under all Higgs fields: it
follows that $2c_1(L_E)$
must be strictly greater than $c_1(\La)$ for $E$ to form a stable pair.

The converse is straightforward:  we suppose that $g=0$, $2c_1(L_E) > c_1(\La)$
and $L_E$ is such
that the bounds in \refle{chi}, \refpa{xpos} are attained and, exactly as in
\refpr{g0}, we set $E=L_E\oplus L_E^*\Lambda$. We write the Higgs field
according to the
decomposition as $\phi =\left( \begin{array}{cc}0 & u\\ v &
-0\end{array}\right)$. Since $g=0$,
the fact that $L_E$ is such that the bounds in \refle{chi}, \refpa{xpos} are
attained means that
$h^0(KL_E^{-2}\La) = n - n_0 -2 \ge 1$ and $h^0(KL_E^{2}  \Lambda^*) = 0$.
Hence $v$ can be chosen
non-zero so that $E$ exists in a stable pair and $u=0$ so that $L_E^*\Lambda$
is invariant by all
$\phi$, as required.

Finally, examples where the bounds in \refle{chi}, \refpa{xpos} are attained
exist by
\refle{bounds}. \epf

\bse{The Yang-Mills-Higgs Equations and Moduli}{ymh}

We now prove an equivalence between stable Higgs $V$-bundles and the
appropriate
analytic objects---irreducible \ymh\ pairs---and use this to give an analytic
construction
of the moduli space. Throughout this section $M$ is an orbifold Riemann surface
of negative Euler
characteristic, equipped with a normalised volume form, $\Omega$, and $E$ is a
smooth
rank 2 $V$-bundle over $M$ with a fixed Hermitian metric.

\bsu{The \ymh\ Equations}{ymhymh}

Given the fixed Hermitian metric on $E$, holomorphic structures correspond to
unitary
connexions. Let $\phi$ be a Higgs field with respect to $A$, \ie a Higgs field
on ${E}_A$ or
satisfying $\o\partial_A\phi =0$.  We call the pair $(A,\phi)$ a \de{Higgs
pair}.  (With
the unitary structure understood Higgs pairs are entirely equivalent to the
corresponding
Higgs $V$-bundles and so we can talk about stable Higgs pairs, isomorphisms of
Higgs pairs
and so on.)  (From some points of view it is more natural to consider the
holomorphic structure as fixed and the unitary structure as varying.  Of course
the two
approaches are equivalent.)

We impose determinant-fixing conditions in what follows; they are not essential
but they remove some redundancies associated with scalar automorphisms (see
\refpr{stable regular}), tensoring by line $V$-bundles and so on.  We have
already made the assumption that the Higgs field $\phi$ fixes determinants in
the sense that it is trace-free; the other determinant-fixing conditions are
defined as follows.  A unitary structure on $E$ induces one on the determinant
line $V$-bundle $\Lambda$.  With this fixed and a choice of isomorphism class
of
holomorphic structure on $\Lambda$, there is a unique (up to unitary gauge)
unitary connexion on $\Lambda$ which is compatible with the class of
holomorphic
structure and is Yang-Mills, \ie has constant central curvature $-2\pi
\ci\,c_1(\Lambda)\Omega$.  Fix one such connexion and denote it $A_\La$.  We
say
that a unitary connexion or holomorphic structure on $E$ has \de{fixed
determinant} if it induces this fixed connexion or holomorphic structure in the
determinant line $V$-bundle.  (On the other hand if we fix the holomorphic
structure then we can choose a \hym\ metric on the determinant line $V$-bundle
and fix the determinant of our metrics by insisting that they induce this
metric.)

Given a unitary connexion $A$ the trace-free part of the curvature is
$F_A^0 =_{\rm def} F_A + \pi \ci\, c_1(\Lambda)\Omega I_E$,
by the Chern-Weil theory.  We say that a Higgs pair $(A,\phi)$ (with fixed
determinants
understood) is \de{Yang-Mills-Higgs} if
\beql{hymh condition}
\begin{array}{rcl} F_A^0 + [\phi,\phi^*] &=& 0 \quad{\rm and}\\
\o\partial_A\phi &=& 0.
\end{array}
\eeql
(For a Hermitian metric varying on a fixed Higgs $V$-bundle this is
the condition for the metric to be \hymh .) The
involution $\phi\mapsto \phi^*$ is a combination of the conjugation $dz\mapsto
d\overline{z}$ and taking the adjoint of an endomorphism with respect to the
metric.  The
second part of the condition merely reiterates the fact that $\phi$ is
holomorphic with
respect to the holomorphic structure induced by $A$.  Of course if $\phi=0$
then
\refeq{hymh condition} is just the Yang-Mills equation (see \cite{ab82,fs92})
and we say
that $A$ is \de{Yang-Mills}.  An existence theorem for Yang-Mills connexions
in
stable $V$-bundles, generalising the Narasimhan-Seshadri theorem from the
smooth
case \cite{do83}, is given in \cite{fs92}.
The first half of our correspondence between stable Higgs $V$-bundles and
\ymh\
pairs is not difficult; again a result of Hitchin \cite[theorem 2.1]{hi87}
generalises easily.
\bpr{stable}
Let $M$ be an orbifold Riemann surface with negative Euler characteristic.
If $(A,\phi)$ is a \ymh\ pair (with respect to the fixed unitary structure on
$E$ and with
fixed determinants) then the pair $(A,\phi)$ is stable unless it has a
$U(1)$-reduction,
in which case it is polystable.
\epr
We call a pair with a $U(1)$-reduction, {\em as a pair}, \de{reducible};
otherwise the
pair is \de{irreducible}. Notice that a reducible pair is \ymh\ if and only if
the
connexions in the two line $V$-bundles are Yang-Mills.

Define the \de{gauge group} ${\cal G}(E)$ to be the group of unitary
automorphisms of $E$ (fixing the base).  This acts on Higgs fields by
conjugation and has a natural action on $\o\partial$-operators such that the
corresponding Chern connexions transform in the standard way.  Thus this action
fixes the determinant line $V$-bundle, acts on the set of Higgs $V$-bundles by
isomorphisms and takes one \ymh\ pair to another.  We also consider the
\de{complexified gauge group} ${\cal G}^c(E)$ of complex-linear automorphisms
of
$E$ (fixing the base).  Again this acts on Higgs $V$-bundles by isomorphisms.

Isomorphic Higgs $V$-bundle structures are precisely those that lie in the same
${\cal G}^c(E)$-orbit.  Notice that \refpr{stable regular} implies that ${\cal
G}^c(E)$ acts freely (modulo scalars) on the set of stable Higgs $V$-bundles.
(If we think of the Higgs $V$-bundle $({E},\phi)$ as fixed and the Hermitian
metric as variable then ${\cal G}^c(E)$ acts transitively on the space of
Hermitian metrics.)  Once again we easily obtain a uniqueness result due to
Hitchin \cite[theorem 2.7]{hi87} in the smooth case.
\bpr{regular pairs}
Let $({E}_1,\phi_1)$ and $({E}_2,\phi_2)$ be isomorphic Higgs $V$-bundles with
fixed
determinants, with Chern connexions $A_1$ and $A_2$ and the same underlying
rank 2
Hermitian $V$-bundle.  Suppose that the Higgs pairs $(A_1,\phi_1)$ and
$(A_2,\phi_2)$ are
both \ymh.  Then $({E}_1,\phi_1)$ and $({E}_2,\phi_2)$ are gauge-equivalent
(\ie there is
an element of ${\cal G}(E)$ taking one to the other).  \epr

\bsu{An Existence Theorem for \ymh\ Pairs}{ymhexi}

A version of the Narasimhan-Seshadri theorem for stable Higgs $V$-bundles
(essentially a converse to \refpr{stable}) can be proved directly for
orbifolds,
extending the arguments of \cite{do83,hi87}.
\bth{Narasimhan-Seshadri}
Let $E\to M$ be a fixed $U(2)$ $V$-bundle over an orbifold Riemann surface of
negative
Euler characteristic.  If $(A,\phi)$ is a polystable Higgs pair with fixed
determinant on
$E$ then there exists an element $g \in {\cal G}^c$ of determinant 1, unique
modulo
elements of $\cal G$ of determinant 1, such that $g(A,\phi)$ is \ymh.
\eth
We shall deduce the theorem from the ordinary case by equivariant arguments in
\refsu{ymhequ}, though there is some advantage to a direct proof, as
an appeal to Fox's theorem is avoided and uniformisation results from the
following corollary, proved as in \cite[corollary 4.23]{hi87}.
\bco{negative curvature}
If $M$ is an orbifold Riemann surface of negative Euler characteristic
then $M$ admits a unique compatible metric of constant sectional
curvature -4.
\eco
\bpf
We define a stable Higgs $V$-bundle by equipping $E=K\oplus 1$
with the Higgs field
\beq
\phi=\left(\begin{array}{cc}
0 & 0 \\
1 & 0
\end{array}\right).
\eeq
We fix a \hym\ metric on $\Lambda^2E$. From \refth{Narasimhan-Seshadri} we
have a \hymh\ metric $h$ on $E$. Exactly as in \cite[corollary 4.23]{hi87},
this must
split and we obtain a metric on $K$ such that the dual metric in the tangent
bundle has constant sectional curvature -4.
\epf

\bsu{The \ymh\ Moduli Space}{ymhmod}

We now construct the moduli space of irreducible \ymh\ pairs, beginning with a
brief discussion of reducible \ymh\ pairs.  Let $(A,\phi)$ be a reducible
\ymh\ pair on $E$.  The reduction means that there is a splitting of $E$ into a
direct sum
$E=L\oplus L^*\La$, where $L$ and $L^*\La$ have the same degree, with
respect to which $A$ and $\phi$ are diagonal---the resulting Higgs $V$-bundle
is
polystable but not stable.  The isotropy group of the pair $(A,\phi)$ is $S^1$
or $SU(2)$
according to whether the two summands are distinct or identical; since $\phi$
is
trace-free the latter is only possible if $\phi=0$.

Let us now consider the question of the existence of reductions.  Obviously the
essential prerequisite is that $L$ exists such that $L$ and $L^*\La$ have the
same degree.  If $a$ denotes the least common multiple of the $\al_i$'s then
the
degrees of line $V$-bundles have the form $s/a$ for $s\in \Z$ and all $s$
occur.
Thus a necessary condition for a reduction is that $c_1(\La) = s/a$ with $s$
even.  However, even when $s$ is even, there is a further constraint:  the
isotropy of $E$ is fixed and, as before, the isotropy of $L$ must be described
by an isotropy vector $(\epsilon_i)$ with $c_1(L) \equiv
\isum\{\epsilon_i(x'_i-x_i)+(x'_i+x_i)\}/2\al_i \pmod{\Z}$ and so the isotropy
may imply a constraint to finding $L$ with appropriate $c_1(L)$.  For general
$M$ and $E$ it is impossible in \lq most' cases (see \cite{fs92} for details).

{}From now on we make the assumption that the isotropy of $M$ and the degree
and
isotropy of $E$ are such that there are no reducible \ymh\ pairs on $E$.

We outline the deformation theory to show that the moduli space is a
finite-dimensional manifold. (For the purposes of this outline we
suppress the use of Sobolev spaces---this is standard; see \eg \cite{pa'82}.)
Fix an irreducible \ymh\ pair $(A,\phi)$. The `deformation
complex' at $(A,\phi)$ is then the following elliptic complex:
\beql{deformation}
0\to \Gamma(\frak{su}(E)) \stackrel{d_1}{\to} \Gamma(\frak{su}^1(E)) \oplus
\Omega^{1,0}(\frak{sl}(E))
\stackrel{d_2}{\to} \Gamma(\frak{su}^2(E)) \oplus \Omega^{1,1}(\frak{sl}(E))
\to
0,
\eeql
where $\frak{su}^k(E)$ denotes the bundle of skew-adjoint $k$-forms with values
in the trace-free endomorphisms of $E$ and $\frak{sl}(E)$ denotes the bundle
of trace-free endomorphisms of $E$.
Here $d_1$, giving the linearisation of
the action, is given by
$$
d_1 : \psi \mapsto (d_A\psi,\ [\phi,\psi]))
$$
and $d_2$, giving the linearisation of the Yang-Mills-Higgs equations, by
$$
d_2 :  (A',\phi') \mapsto (d_AA' + [\phi',\phi^*] +[\phi,\phi'^*],\
\o\partial_A\phi' + [(A')^{0,1},\phi]).
$$

We use the orbifold Atiyah-Singer index theorem \cite{ka81} to calculate the
index of \refeq{deformation} as $6(g-1) +
2(n-n_0)$. We note that the zeroeth and second cohomology groups, $H^{0}$ and
$H^{2}$, of the complex vanish---for $H^{0}$ this follows from the
irreducibility of $(A,\phi)$ and for $H^{2}$ the duality argument given by
Hitchin will suffice.  Hence the first cohomology group has dimension $ 6(g-1)
+
2(n-n_0)$.  Moreover the Kuranishi method shows that a neighbourhood of zero in
$H^1$ is a local model for the moduli space and hence
the moduli space is a smooth complex manifold of dimension $ 6(g-1) +
2(n-n_0)$.
\bth{moduli} Let $M$ be an orbifold Riemann surface of negative Euler
characteristic and $E\to M$ a fixed complex rank 2 $V$-bundle.
\begin{enumerate}
\item Suppose that $E$ is equipped with a Hermitian metric and admits no
reducible \ymh\ pairs.  Then the moduli space of \ymh\ pairs on $E$ with fixed
determinants, ${\cal M}(E,A_\La)$, is a complex manifold of dimension $
6(g-1) + 2(n-n_0)$.
\item Suppose that $E$ admits no Higgs $V$-bundle structures
which are polystable but not stable.  Then the moduli space of stable Higgs
$V$-bundle
structures on $E$ with fixed determinants is a complex manifold of
dimension $6(g-1) + 2(n-n_0)$.
\end{enumerate}
\eth
\bre{roots}
In the smooth case there are essentially only two moduli spaces (of which only
one is smooth), according to the parity of the degree.  In the orbifold case,
how many moduli spaces are there?  Clearly it is sufficient to consider only
one
topological $\Lambda$ in each class under the equivalence $\Lambda \sim \Lambda
L^2$, for any topological line $V$-bundle $L$---`square-free' representatives
for each class will be discussed in \refsu{reprep}.  A further subtlety in the
orbifold case is the possibility of non-trivial topological square roots of the
trivial line $V$-bundle, or simply \de{topological roots}:  if $L$ is a
topological root then there is a map on moduli $\cm(E,A_\Lambda)
\leftrightarrow
\cm(E\otimes L,A_\Lambda)$ by tensoring by $L$, which fixes $\Lambda$ but
alters
the topology of $E$.  For $L$ to be a topological root necessarily $c_1(L)=0$
and $L$ has `half-trivial' isotropy, \ie the isotropy is 0 or $\al/2$ at each
marked point.  If we consider topological line $V$-bundles of the form $L=
\otimes_{\al_i {\rm\ even}}L_i^{\delta_i\al_i/2}$, for $\delta_i \in
\Z$ where the $L_i$ are the point $V$-bundles of \refsu{orbdiv}, then it
is clear that $L$ is a topological root provided $c_1(L)=\sum\delta_i/2 =0$.
If
we let $n_2$ denote the number of marked points where the isotropy is even,
then, provided $n_2 \ge 1$, there are $2^{n_2-1}$ topological roots.  It
follows
that for each topological $\Lambda$, if $n_2 \ge 1$, there will be $2^{n_2-1}$
different topological $E$'s giving essentially the same moduli space.  We will
see another manifestation of this in \refsu{reprep}.  \ere

Recall that the tangent space to the moduli space is given by the first
cohomology of the
deformation complex \refeq{deformation}, \ie by $\ker{(d_1^*)}\cap\ker{(d_2)}$.
This space
admits a natural $L^2$ metric and, just as in
\cite[theorems 6.1 \& 6.7]{hi87}, we have the following result.
\bpr{metric}
Let $E$ be a fixed rank 2 Hermitian $V$-bundle over an orbifold
Riemann surface of negative Euler characteristic and suppose that $E$ admits no
reducible
\ymh\ pairs.  Then the natural $L^2$ metric on the moduli space ${\cal
M}(E,A_\La)$ is
complete and hyper-K\"ahler.
\epr

\bsu{The \ymh\ equations and Equivariance}{ymhequ}

Here we sketch how many {\em but not all} of the previous results of this
section can be treated by equivariant arguments.  Further details for this
subsection can be found in \cite{na91}.

An orbifold Riemann surface with negative Euler characteristic, $M$, has a
topological orbifold covering by a surface \cite{sc'83} and so its universal
covering is necessarily a surface with negative Euler characteristic.
Pulling-back the complex structure we find that the universal covering is
necessarily $D^2$, the unit disk, with $\pi_1^V(M)$ a group of automorphisms
acting properly discontinuously.  In other words $\pi_1^V(M)$ is a co-compact
Fuchsian group or, in the terminology of \cite{fo52}, an \de{$F$-group}.

Thinking of $D^2$ as the hyperbolic upper half-plane or Poincar\'e disk, the
elements of $\pi_1^V(M)$ act by orientation-preserving isometries and so
we get a compatible Riemannian metric of constant sectional curvature on $M$.
This is just \refco{negative curvature}.  In this context we need the following
result of \cite{fo52}.
\bprn{Fox}{Fox}
If $\G$ is an $F$-group then $\G$ has a normal subgroup of finite
index, containing no elements of finite order.
\epr
\bco{smooth covering}
Let $M$ be an orbifold Riemann surface with negative Euler characteristic. Then
there
exists a smooth Riemann surface, $\lift{M}$, with negative Euler
characteristic, together
with a
finite group, $F$, of automorphisms of $\lift M$, such that $M=F\backslash\lift
M$.
\eco
The important point here is that the covering is {\em finite} and hence
$\lift M$ is compact.

The existence result of \refth{Narasimhan-Seshadri} follows from the
corresponding result on $\lift M$, \cite[theorem 4.3]{hi87}, using an averaging
argument (compare \cite{gp91}).  We will always use the notation that objects
on
$\lift M$ pulled-back from $M$ under the covering map $\lift M\to M$ will be
denoted by a `hat'; $\lift{\ \ }$.  In this notation the pull-back of a
$V$-bundle $E \to M$ becomes $\lift E \to \lift M$, and so on.  For the
equivariant argument it is easiest to fix the Higgs $V$-bundle structure on $E$
and vary the metric; therefore, rather than suppose that a {\em Hermitian}
structure on $E$ is given, we temporarily suppose that a {\em holomorphic}
structure on $E$ (and hence on $\lift E$) is given.  We will show that if
$(E,\phi)$ is stable then $(\lift E,\lift\phi)$ is polystable and admits a
\hymh\ metric which is $F$-invariant and so descends to the required metric on
$E$.
\bpr{polystable}
Let $(E,\phi)$ be a stable Higgs $V$-bundle and let $(\lift E,\lift\phi)$ be
the
pull-back to $\lift M$. Then $(\lift E,\lift\phi)$ is polystable.
\epr
\bpf
Suppose first that $(\lift E,\lift\phi)$ is {\em not semi-stable}.  Then there
is a unique
destabilising Higgs sub-$V$-bundle $L=L_{\lift E}$ and the action of $F$ cannot
fix $L$.
Therefore for some $f\in F$ we have that $f(L) \ne L$.  However $f(L)$ is a
Higgs
sub-$V$-bundle of $(\lift E,\lift\phi)$ (because $\lift\phi$ commutes with the
action of $f\in F$) and has the same degree as $L$.  This contradicts the
uniqueness of
$L$. So $(\lift E,\lift\phi)$ is semi-stable.  Suppose it is {\em not stable}.
Then again there
is a destabilising Higgs sub-$V$-bundle $L= L_{\lift E}$ (not necessarily
unique).
As before $L$ cannot be fixed by $F$ and so we obtain, for some $f\in F$, a
Higgs
sub-$V$-bundle $f(L) \ne L$ of the same degree as $L$. Let $g:f(L)\to \lift
E/L$
be the composition of the inclusion of $f(L)$ into
$\lift E$ with the projection onto $\lift E/L$: $g$ is a homomorphism
between two line bundles of the same degree and hence either
zero or constant.  Since $f(L) \ne L$ the map $g$ cannot be zero and hence
$f(L)=\lift
E/L$.  Since $f(L)$ is actually a Higgs sub-$V$-bundle, $(\lift E,\lift\phi)$
is
a direct sum $(L \oplus f(L),\lift\phi_{L}\oplus\lift\phi_{f(L)})$ and so
is polystable as claimed.  \epf

\bpr{metrics exist}
Let $(E,\phi)$ be a stable Higgs $V$-bundle and let $(\lift E,\lift\phi)$ be
the pull-back to $\lift M$. Then the  polystable Higgs $V$-bundle $(\lift
E,\lift\phi)$
admits a \hymh\ metric which is $F$-invariant (and unique up to scale).
\epr
\bpf
Certainly $(\lift E,\lift\phi)$ admits a \hymh\ metric (by \refpr{polystable}
and \cite[theorem 4.3]{hi87}).  By averaging, the \hymh\ metric can be supposed
$F$-invariant.  \epf
An $F$-invariant \hymh\ metric descends to $(E,\phi)$, where it
trivially still satisfies the \hymh\ condition. We can satisfy the
determinant-fixing condition by a choice of scalar multiple and so we
obtain the desired existence result---\refth{Narasimhan-Seshadri}.

Suppose again that a Hermitian, rather than holomorphic, structure on $E$ is
given.  We recall that Hitchin proves that if $\lift E$ has odd degree then
there is a smooth moduli space ${\cal M}(\lift E,\lift A_\La)$ of complex
dimension $6(\lift g-1)$.  The pull-back map $(A,\phi) \mapsto (\lift
A,\lift\phi)$ defines a map from Higgs pairs on $E$ to $F$-invariant Higgs
pairs
on $\lift E$---what can be said about the corresponding map on moduli?  Suppose
that $(A,\phi)$ is an irreducible \ymh\ pair on $E$.  The first point to note
is
that $(\lift A,\lift\phi)$ may be reducible, by the analogue of
\refpr{polystable} for pairs.  For simplicity, we will ignore this possibility
in our discussion---we suppose that there are topological obstructions to the
existence of reducible \ymh\ pairs on $\lift E$.

\ble{regular lifts} Suppose that $(A,\phi)$ is an irreducible \ymh\ pair
on $E$ with an irreducible lift. Suppose further that for some
$g\in \lift{\cal G}$, of determinant 1, $g(\lift A,\lift\phi)$ is
$F$-invariant. Then
$f^{-1}gf=\pm g$ for all $f\in F$.
Conversely, given $g\in \lift{\cal G}$ of determinant 1 such that $f^{-1}gf=\pm
g$ for all
$f\in F$, $g(\lift A,\lift\phi)$ is irreducible and $F$-invariant.
\ele
\bpf
Since $(\lift A,\lift\phi)$ is $F$-invariant we
know that $fd_{\lift A} = d_{\lift A}f \and f\lift\phi = \lift\phi f $ for any
$f\in F$. Since the same is also true of $g(\lift A,\lift\phi)$ It follows
that
$d_{\lift A} = (g^{-1}f^{-1}gf)(d_{\lift A})(f^{-1}g^{-1}fg)$ and similarly
for the Higgs field. Since $(A,\phi)$ is a stable pair it follows
(\refpr{stable
regular}) that $\pm g = f^{-1}gf$. The converse is clear.\epf

Let $\lift\cg^F$ be the subgroup of $\lift\cg$ consisting of $F$-invariant
elements of determinant 1 and let $\lift{\cal G}^{\pm}$ denote that of elements
$g\in \lift{\cal G}$ of determinant 1 such that, for all $f\in F$,
$f^{-1}gf=\pm
g$.  Clearly either $\lift{\cal G}^{\pm}=\lift\cg^F$ or $\lift\cg^F <
\lift{\cal
G}^{\pm}$ with even index.  (In fact these groups will be equal under quite
mild
hypotheses, which amount to the vanishing of a certain equivariant
$\Z_2$-characteristic class---see \cite{na91} and compare \cite[proposition
1.8,
part iii)]{fs92}.)  If these groups are unequal then $f^{-1}gf = -g$ for some
$f\in F$ and $g\in \lift{\cal G}$ of determinant 1---but such a $g$ cannot be
close to $\pm 1$ and so does not enter the local description of the moduli
space
(compare \cite[theorem 4.1]{pa'82}). At an irreducible $F$-invariant pair
$(\lift A,\lift\phi)$ the group $F$ acts on the
deformation complex. The pull-back map induces a commutative diagram of
deformation complexes and it follows immediately that ${\cal M}(E,A_\La)$
covers a submanifold of ${\cal M}(\lift E,\lift A_\La)$ with covering group
$\lift{\cal G}^{\pm}/\lift{\cal G}^{F}$.
\bth{sub}
Let $M$ be an orbifold Riemann surface of negative Euler
characteristic and $E\to M$ a fixed complex rank 2 $V$-bundle. Let $\lift E$
be
the pull-back of $E$ under the identification $M = F\backslash \lift M$
of \refco{smooth covering}.
\begin{enumerate}
\item Suppose that $E$ is equipped with a Hermitian metric and $\lift E$ with
the pulled-back metric and that $E$ admits no
reducible \ymh\ pairs. If $\lift E$ has odd degree then, under pull-back,
the moduli space of \ymh\ pairs with fixed determinants on $E$, ${\cal
M}(E,A_\Lambda)$, covers a
submanifold of the corresponding moduli space on $\lift E$ with
covering group $\lift{\cal G}^{\pm}/\lift{\cal G}^{F}$ (with
$\lift\cg^\pm$ and $\lift\cg^F$ as above).  If $\lift E$ has even degree
then this remains true for those classes of Higgs pairs which are irreducible
on $\lift
E$.
\item Suppose that $E$ admits no Higgs $V$-bundle structures
which are polystable but not stable.  If $\lift
E$ has odd degree then, under pull-back, the moduli space of stable Higgs
$V$-bundle
structures with fixed determinants on $E$ covers a submanifold of the
corresponding moduli
space on $\lift E$ with covering group $\lift{\cal G}^{\pm}/\lift{\cal G}^{F}$
(with
$\lift\cg^\pm$ and $\lift\cg^F$ as above). If
$\lift E$ has even degree then this remains true for those classes of Higgs
$V$-bundle
structure which are stable on $\lift E$.
\end{enumerate}
\eth

Notice that in the case when $\lift M$ is a hyperelliptic surface of
genus 2 branched over 6 points of the Riemann sphere then the dimensions of the
two moduli spaces are equal (a simple arithmetic check shows that this is the
only case where this happens).

\bse{The topology of the moduli space}{top}

We now give some results on the topology of the moduli space using the Morse
function $(A,\phi)\stackrel{\mu}{\to}||\phi||_{L^2}^2$, following \cite[\S
7]{hi87}.  Notation and assumptions remain as before; in particular, we suppose
that $E$ admits no reducible \ymh\ pairs, so that the moduli space $\cm = \cm
(E,A_\Lambda)$ is smooth and recall the definitions of the integers $n_\pm$
and
$l$ from \refsu{orbint}.

The function $(A,\phi)\stackrel{\mu}{\to}||\phi||_{L^2}^2=2\ci \int \tr(\phi
\phi^*)$ is invariant with respect to the circle action
$e^{\ci\theta}(A,\phi)=(A,e^{\ci\theta}\phi)$ and $d\mu(Y) = -2\ci
\omega_1(X,Y)$ where $X$ generates the $S^1$-action and $\omega_1$ is as in
\cite[\S 6]{hi87}.  The map $\mu$ is proper and there's an extension of
\cite[proposition 7.1]{hi87}.  To describe it we need to consider pairs
$(m,(\epsilon_i))$ where $m$ is an integer and $(\epsilon_i)$ is an isotropy
vector---such pairs describe topological sub-$V$-bundles of $E$, with isotropy
described by $(\epsilon_i)$ and degree $m + \isum
\{\epsilon_i(x'_i-x_i)+(x'_i+x_i)\}/(2\al_i)$ (see \eg \refre{subbundles}).

\bth{Morse}  Let $E$ be a fixed rank 2 Hermitian $V$-bundle over an orbifold
Riemann surface of negative Euler characteristic and suppose that $E$ admits no
reducible
\ymh\ pairs.  If $g=0$ then
suppose that
$n-n_0\ge 3$.  Let $\mu$ be as above: then, with the notations established
above,
\begin{enumerate}
\item\label{critical values} $\mu$ has critical values 0 and
$2\pi\{ 2 m -l + \isum \{\epsilon_i(x'_i-x_i)/\al_i\} \}$ for an integer $m$
and
isotropy vector $(\epsilon_i)$ with
\beq
l < 2m + \isum \frac{\epsilon_i(x'_i-x_i)}{\al_i} \le  l + 2g - 2 +
\isum\frac{\epsilon_i(x_i' - x_i)}{\al_i} + n_-;
\eeq
\item the minimum $\mu^{-1}(0)$ is a non-degenerate critical manifold of index
0 and is
diffeomorphic to the space of stable $V$-bundles with fixed determinants and
\item the other critical manifolds are also
non-degenerate and are $2^{2g}$-fold coverings of
$S^r\wo M$, where $r = l- 2m +2g -2 + n_-$.  Moreover, they
are of index $2\{2m -l + g - 1 + n_+ \}$.
\end{enumerate}
\eth
\bpf
The critical points are the fixed points of the induced circle action on $\cm$.
Because we are taking quotients by the gauge group, these correspond to pairs
$((A,\phi),\la)$ where $\la :  S^1 \to {\cal G}$ such that, for all $\theta$,
$\la(e^{\ci \theta})d_A\la(e^{-\ci \theta}) = d_A$ and $\la(e^{\ci
\theta})\phi\la(e^{-\ci \theta}) = e^{\ci \theta}\phi$.  If $\phi=0$ then,
holomorphically, we simply get stable $V$-bundles.  If $\phi\ne 0$ then
certainly $\la(e^{\ci \theta})\ne 1$ for $\theta \not\equiv 0 \pmod{2\pi}$.
The
first equation now implies that the stabiliser $\cg_A$ is non-trivial and $A$
is
reducible to a $U(1)$-connexion.  Consequently, as a holomorphic $V$-bundle,
$E$
is decomposable (so, in particular, not stable) and can be written $L\oplus L^*
\Lambda$.  If we write $\phi = \left( \begin{array}{cc} t & u\\v & -t
\end{array}\right)$ and $\la(e^{\ci \theta})= \left( \begin{array}{cc}
\mu_\theta & 0\\0 & \mu_\theta^{-1} \end{array}\right)$ with respect to this
splitting then the second equation implies $t=0$ and either $u=0$ or $v=0$.
Replacing $L$ by $L^* \Lambda$ if necessary, we can suppose that $u=0$ and that
$v\in H^0(KL^{-2} \Lambda)$---$v$ is holomorphic from the self-duality
equations.

The remaining term of the \ymh\  equations is
$*(F_A + [\phi,\phi^*] )= -\pi\ci d I_E$.  Writing $*F_{A_L} = *F-\pi\ci d$,
in
terms of the above decomposition, so that $*F_{A_{L^*\La}} =
-*F-\pi\ci d$, we find that $F = v \wedge \o v$ and
\beq
\deg L
=\frac{\ci}{2\pi}\int(F - *\pi\ci c_1(\Lambda)) = \frac{\ci}{2\pi}\int (v\wedge
\o v) +
\frac{c_1(\Lambda)}2 = \frac{\mu}{4\pi} + \frac{c_1(\Lambda)}{2}.
\eeq
Since $\mu> 0$ for $\phi\ne 0$, we have $2\deg L > c_1(\Lambda)$ and $L=L_E$,
the
destabilising
sub-$V$-bundle of
$E$.  Moreover, because $v\ne 0$ we must have $h^0(KL^{-2}\Lambda)\ge
1$ (compare \refth{stable pairs}).

Now, for any $(m,(\epsilon_i))$ let $L_{(m,(\epsilon_i))}$ be the
corresponding
topological sub-$V$-bundle of $E$. Consider pairs $(m,(\epsilon_i))$ with
$2c_1(L_{(m,(\epsilon_i))})> c_1(\Lambda)$ and
set $L=L_{(m,(\epsilon_i))}$ and $E=L\oplus L^*\La$.
This occurs as a stable pair $(E,\phi)$ provided $L$ admits a holomorphic
structure with $h^0(KL^{-2}\Lambda) \ge 1$, and the Higgs field $\phi$ is then
given by $v\in H^0(KL^{-2}\Lambda)\setminus\{0\}$ (compare \refth{stable
pairs},
\refpa{nsde} and \refpr{counter}).

To see whether a given topological $L=L_{(m,(\epsilon_i))}$ admits an
appropriate holomorphic
structure we use our results from \refsu{higalg}: by \refle{chi} we have
$\chi(KL^{-2}  \Lambda) =
l -  2m + g - 1  + n_-$. It follows that $r=
c_1 (\wo{ KL^{-2}\La }) = l -  2m + 2g -2  + n_-$.
Hence, supposing that $r  \ge 0$, for each effective
(integral) divisor of divisor order $r$ (if $r=0$ then for the empty divisor)
we obtain a
holomorphic structure on $\wo{K L^{-2}\La}$ with a holomorphic section
determining the divisor (determined up to multiplication by elements of
$\C^*$).  Hence we get a
holomorphic structure on $K L^{-2}\La$ with holomorphic section $v$ and all
holomorphic sections arise in this way.  Placing a corresponding holomorphic
structure on
$L$ requires a choice of holomorphic square root and there are $2^{2g}$ such
choices.  For each root $L$ the pair $(E,v) =(L\oplus L^*\La,v)$ is clearly
stable by construction. The section $v$ is determined by the divisor up to a
multiplicative constant
$\lambda\ne 0$ but $(L\oplus L^*\La,v)$ and $(L\oplus L^*\La,\lambda v)$ are in
the same orbit under
the action of the complexified gauge group and hence
equivalent.  Two distinct divisors determine distinct stable pairs so that we
have the critical set
is a $2^{2g}$-fold covering of the set of effective divisors of degree $r =
l -  2m + 2g -2  + n_-$; that is, a $2^{2g}$-fold covering of
$S^r\wo M$ (a point if $r=0$).

Let $E=L\oplus L^*\La$ for $L=L_{(m,(\epsilon_i))}$, as above.  The subset
$U=\left\{ \phi \in
H^0(\End_0(E)\otimes K)\ :\ (E,\phi)\mbox{ is stable } \right\}$ is acted upon
freely by
$\Aut_0(E)/\{\pm 1\}$, where $\Aut_0(E)$ are the holomorphic automorphisms of
determinant 1 (see
\refpr{stable regular}).  The quotient $U/(\Aut_0(E)/\{\pm 1\})$ is a complex
manifold of dimension
$3g-3+n-n_0$.  So through each point $P\in \cm$ there passes a
$(3g-3+n-n_0)$-dimensional isotropic
complex submanifold $U/(\Aut_0(E)/\{\pm 1\})$, invariant under $S^1$:  it is
thus {\em Lagrangian}.
Suppose $P\in \cm$ is fixed under the $S^1$-action and $P=(E,\phi)$, where
$E=L\oplus L^*\La$, $\phi
= \left( \begin{array}{cc}0 & 0 \\ v & 0 \end{array}\right)$, as above.  The
homomorphism $\lambda$
is given by $\lambda(\theta) = \left( \begin{array}{cc}e^{-\ci\theta/2} & 0 \\
0 & e^{\ci\theta/2}
\end{array}\right)$ with respect to this decomposition.  Now $\End_0(E) =
L^{-2}\La\oplus
L^2\La^*\oplus \C$ and $\lambda(\theta)$ acts as
$(e^{\ci\theta},e^{-\ci\theta},1)$.  Hence
$\lambda(\theta)$ acts with negative weight solely on $H^0(KL^2\La^*) \subset
H^0(\End_0(E)\otimes
K)$.  As $\lambda(\theta)$ acts on $\phi$ by multiplication by $e^{\ci\theta}$
there are no negative
weights on $H^0(\End_0(E)).\phi$ and hence we find, as in \cite{hi87}, that the
index is $2
h^0(KL^2\La^*) = 2\{2m -l + g - 1 + n_+ \}$, by
\refle{chi}. \epf

{}From this, the work of \cite{fr59} and general Morse-Bott theory
\cite{ab'99} we can, in
principle, calculate the Betti numbers---see \cite{by}.
We content ourselves with \refco{topology}, below, for which we need the
following preliminary lemma.
\ble{index0}
There is exactly one critical manifold of index 0 and this
is connected and simply-connected.
\ele
\bpf
\refth{Morse} shows that if $g>0$ then the space of stable $V$-bundles
is the only index 0 critical manifold and this is connected and
simply-connected (even when $g=0$) by \cite[theorem 7.11]{fs92}.
When $g=0$, critical manifolds of index 0 other than the moduli of
stable $V$-bundles may occur: these have the form $S^{r}\tilde M \cong
\C\P^{r}$ and so are also connected and simply connected.

It remains to show that exactly one of the possibilities is
non-empty in each case. Making allowances for differences in notation,
the following is implicit in \cite[theorem 4.7]{fs92}:
the space of stable $V$-bundles is empty if and only if there exists a
vector $(\epsilon_i)$ with $n_+ +l \equiv 1 (2)$ and
\beql{empty}
n_+ - \isum\frac{\epsilon_i(x'_i-x_i)}{\al_i} < 1-g.
\eeql
Since the left-hand side of \refeq{empty} is clearly not less than zero we
see that the space of stable $V$-bundles is non-empty whenever $g>0$.

When $g=0$, \refth{Morse} shows that the critical manifolds of index 0 other
than the moduli of stable $V$-bundles consist precisely of the
$V$-bundles considered in \refpr{counter}. The number of such critical
manifolds
is the number of topological types
$L_{(m,(\epsilon_i))}$ satisfying the criteria of \refpr{counter}, which,
using the ideas of \refle{bounds}, is
\beql{count}
\# \left\{ (\epsilon_i)\  :\  n_+ + l\equiv 1 (2) \quad\mbox{and}\quad   n_+
-
\isum\frac{\epsilon_i(x'_i-x_i)}{\al_i}  < 1 \right\},
\eeql
where $(\epsilon_i)$ varies over all isotropy vectors.
Comparing \refeq{count} to \refeq{empty} we see that exactly one
of the two types of critical manifold must occur.
Moreover, we claim that the number in \refeq{count} is at most 1---this
is sufficient to establish the lemma.

To prove the claim suppose, without loss of generality, that $n_0=0$.
Observe that it is an easy exercise to show that
if $t_1, \dots t_{n}\in (0,1)$ are such that
$\sum_{i=0}^{n} t_i <1$ then at most one $t_i$ can be replaced by
$1-t_i$ with the sum remaining less than 1. Let
\beq
t_i = \frac{1+\epsilon_i}{2} - \frac{\epsilon_i(x'_i - x_i)}{\alpha_i}
\eeq
so that $\sum_{i=0}^n t_i = n_+ - \sum_{i=0}^n \epsilon_i(x'_i - x_i)/
\alpha_i$ and changing the sign of $\epsilon_i$ simply sends $t_i$ to
$1-t_i$. The observation applies to show that this sum can be less than
1 for at most two vectors $(\epsilon_i)$ and these cannot have $n_+$ of
the same parity. Hence the count in \refeq{count} is at most 1, as
claimed.
\epf

\bco{topology}
The moduli space $\cm$ is non-compact---except in the case $g=0$ and
$n-n_0=3$ when it is a point---and connected and simply-connected.
\eco
\bpf
The non-compactness follows from the fact that the critical manifolds cannot be
maxima except if $g=0$ and $n-n_0 = 3$.  This is because the critical manifolds
have index $ i = 2\left\{2m -l + g - 1 + n_+ \right\}$ and (real) dimension $2r
= 2\left\{ l - 2m + 2g -2 + n_- \right\}$ and $2r+i = 6g - 6 +2(n-n_0)$, which
is exactly half the (real) dimension of the moduli space.  The connectedness
and
simple-connectedness follow from the analogous facts for the unique critical
manifold of index 0 (\refle{index0}) and the fact that the other Morse indices
are all even and strictly positive.
\epf

\bse{The Determinant Map}{det}

Recall that $M$ is an orbifold Riemann surface with negative Euler
characteristic, with $E\to M$ a fixed $U(2)$-$V$-bundle.  We assume that $E$
admits no reducible \ymh\ pairs so that the moduli space is smooth.

Thinking of the moduli space as a space of stable Higgs $V$-bundles, there is a
holomorphic gauge-invariant map $ (A,\phi) \mapsto \det(\phi) $ which descends
to a holomorphic map $ \det :  \cm(E,A_\Lambda) \to H^0(K^2).$ Hitchin showed
that in the smooth case this map is proper, surjective and makes $\cm$ a
completely integrable Hamiltonian system.  Moreover he showed that when $q \in
H^0(K^2)$ has simple zeros the fibre $\det^{-1}(q)$ is biholomorphic to the
Prym
variety of the double covering determined by $\sqrt{-q}$ \cite[theorem
8.1]{hi87}.  We will see that things are similar but a little more involved in
the orbifold case:  the first significant observation is that $h^0(K^2) = 3g
- 3 + n$---this is half the dimension of the moduli space exactly when $n_0 =
0$.  For this reason it will be useful to suppose that $n_0 = 0$.  (In
\refsu{detred} we will show that the image of the determinant map is contained
in a canonical $(3g - 3 +n-n_0)$-dimensional subspace of $H^0(K^2)$ and thus
all
cases can be reduced to the case $n_0 = 0$.)  In addition, there are two
special
cases which we exclude:  when $g=0$, $n=3$ the determinant map is identically
zero, and when $g=1$, $n=1$ we have a special case which leads to a breakdown
in
our methods---this case is dealt with separately in \refsu{detspe}.

We summarise our results in the following theorem (proofs are for the most part
discussed in the remainder of this section; the details which have been omitted
are exactly as in \cite[\S 8]{hi87}).  We believe that a similar result was
obtained by Peter Scheinost.
\bth{determinant map}
Let $E$ be a fixed
rank 2 Hermitian $V$-bundle over an orbifold Riemann surface of negative Euler
characteristic, with
$n-n_0>3$ if $g=0$.  Suppose further that $E$ admits no reducible \ymh\ pairs.
Then the determinant
map on the moduli space of \ymh\ pairs on $E$ with fixed determinants
\beq
\det :  \cm(E,A_\Lambda) \to H^0(K^2)
\eeq
has the following properties:  \begin{enumerate} \item $\det$ is proper; \item
the image of $\det$
lies in a
canonical $(3g - 3 +n-n_0)$-dimensional subspace $H^0(\b M;K_{\b M}^2)\subseteq
H^0(K^2)$ and
$\det$ surjects onto $H^0(\b M;K_{\b M}^2)$;
\item with respect to $\det :  \cm(E,A_\Lambda) \to H^0(\b
M;K_{\b M}^2)$, $\cm(E,A_\Lambda)$ is a completely integrable Hamiltonian
system; \item for a generic $q$ in
the image of $\det$, the fibre $\det^{-1}(q)$ is biholomorphic to a torus of
dimension
$3g-3+n-n_0$---this can be identified with the Prym variety of the covering
determined by
$q$ except when $g=n-n_0=1$, when it is identified with the Jacobian; \item
$\cm(E,A_\Lambda)$ is a fibrewise compactification of $T^*\cn(E,A_\Lambda)$
with respect
to the map $\det :  T^*\cn(E,A_\Lambda) \to H^0(\b M;K_{\b M}^2)$, where
$\cn(E,A_\Lambda)$ is the
moduli space of Yang-Mills connexions on $E$ with fixed determinants.
\end{enumerate}
\eth

It seems possible to obtain results arguing using orbifold methods but it is
often simpler to
translate this orbifold problem into one about parabolic bundles; we review the
necessary results in
the next subsection.

\bsu{Parabolic Higgs bundles}{parhig}

Recall the basic facts concerning the correspondence between $V$-bundles over
$M$ and parabolic
bundles over $\wo M$ \cite{fs92}.  Let $\wo E$ be a rank $2$ holomorphic vector
bundle over $\wo M$.  A
\de{quasi-parabolic structure} on $\wo E$ is, for each marked point $p \in \{
p_1,\dots,p_n\}$, a flag
in $\wo E_p$ of the form \beq
\wo E_{p} = \C^2 \supset \C \supset 0, &\mbox{ or }& \wo E_p = \C^2
\supset 0.
\eeq
A flag of the second form is said to be \de{degenerate}.  A quasi-parabolic
bundle $\wo E$ is a \de{parabolic bundle} if to each flag of the first form
there is attached a pair of weights, $0\le \la < \la' < 1$ and to each of the
second form there is a single (multiplicity 2) weight $0\le \la = \la' < 1$.
There is a notion of parabolic degree involving the degree of $\wo E$ and the
weights.  A basis $\{ e,e' \}$ for the fibre at a parabolic point is said to
\de{respect the quasi-parabolic structure} if either the flag is degenerate or
$e'$ spans the intermediate subspace in the flag.  An endomorphism of a
parabolic bundle $\psi$ is a \de{parabolic endomorphism} if for each $p$, with
respect to a basis which respects the quasi-parabolic structure, $\psi_p$
satisfies $(\psi_p)_{12} = 0$ whenever $\la < \la'$.

Let $E$ be a rank $2$ holomorphic $V$-bundle over $M$.  Recall
that by convention $x \le x'$ (if we assume that $n_0 = 0 $ then there is
strict inequality).  For a
line $V$-bundle $L$, we can consider the passage $L\mapsto \wo L$
(\refsu{orbdiv}) as a
smoothing process and the construction of parabolic bundles follows similar
lines:  for a marked
point $p$ we consider \beq (E|_{M \setminus\{ p \}}) \cup_\Psi D^2 \times \C^2,
\eeq with clutching
function $\Psi$ given, in local coordinates, by its $\Z_\al$-equivariant
lifting \beql{patch}
\begin{array}{rcl} \lift\Psi :  (D^2 \setminus\{0\}) \times \C^2 &\to& D^2
\times \C^2\\
(z,(z_1,z_2)) &\mapsto& (z^\al,(z^{-x}z_1,z^{-x'}z_2)).  \end{array} \eeql Now
a holomorphic
section of $(D^2 \times \C^2)/(\sigma \times \tau)$ is given by holomorphic
maps $s_j :  D^2\to \C$,
for $j=1,2$, invariant under the action of $\Z_\al$.  As with
\refeq{Taylor}, Taylor's theorem implies that $s_j(z) = z^{x_j}
\wo{s}_j(z^\al)$, where $\wo{s}_j$ is a
holomorphic function $D^2 \to \C$ and we use the temporary notations $x_1=x$
and $x_2=x'$.
Under the map $\Psi$ defined by \refeq{patch} we simply get
a section of $(D^2 \setminus \{ 0 \})\times \C^2$ which is given by the
functions $\wo{s}_j(w)$ and
hence extends to a holomorphic section of $D^2\times \C^2$.  In other words the
map $\Psi$ is an
isomorphism between the sheaves of germs of holomorphic sections.  Repeating
this construction about
each marked point, we get a holomorphic bundle $\wo E\to \wo M$ corresponding
to the holomorphic
$V$-bundle $E\to M$.

In fact $\wo E$ has a natural parabolic structure as follows:  working in our
local coordinates
about a particular marked point (which respect the $V$-structure) we define
weights $\la=x/\al$ and
$\la'=x'/\al$.  Define a flag in $\C^2$ so that the smallest proper flag space
is
the subspace of $\C^2$ on which $\tau$ acts like $\sigma^{x'}$.  The
corresponding quasi-parabolic
structure on $\wo E_p$ is then given by the image of this flag---notice that
this is degenerate if
and only if $x = x'$.  With the weights $\la,\la'$ it is clear that $\wo E$ is
a parabolic bundle.
(Whilst it is not true in general that $\La^2\wo E = \wo{\Lambda}$, the
bundle $\La^2\wo E$ is determined by $\Lambda$ and the isotropy so that our
determinant-fixing condition on $E$ translates to one on $\wo E$.)
We quote the following result of \cite{fs92}.
\bprn{Furuta-Steer}{V-parabolic}
For a fixed
orbifold Riemann surface $M$, the correspondence $E \mapsto \wo E$ gives a
bijection between
isomorphism classes of rank 2 holomorphic $V$-bundles and those of rank 2
parabolic bundles over
$\wo M$ with rational weights of the form $x/\al$.  Moreover, the induced map
$\co(E) \mapsto
\co(\wo E)$ is an isomorphism of analytic sheaves.  \epr

Now consider what happens to Higgs fields under the passage $E \mapsto \wo E$:
we use a local uniformising coordinate $z$, centred on a given marked point,
and let $w=z^\al$ be the
local holomorphic coordinate on $\wo M$.  There is a Taylor series expansion
as before:  if $\phi$ is a Higgs field on $E$ then in our local coordinates
\beql{Taylor Higgs}
\phi_{ij}dz &=&\left\{\begin{array}{ll}
z^{x_i-x_j-1}\wo\phi_{ij}(z^\al)dz&\quad\mbox{ if }x_i>x_j\and\\ z^{\al +
x_i-x_j-1}\wo\phi_{ij}(z^\al)dz&\quad\mbox{ if }x_i \le x_j,\\
\end{array}\right.  \eeql
where $\wo\phi_{ij}$ are holomorphic functions and we again use the temporary
notations
$x_1=x$ and $x_2=x'$.

To transfer this across to $\wo E$ simply notice that away from the marked
point
the clutching
function $\Psi$ defined by \refeq{patch} is a bundle isomorphism and so acts on
the Higgs field by
conjugation.  Conjugating by $\Psi$ we obtain \beql{Taylor Higgs 2}
\phi_{ij}^\Psi dz &=&
z^{x_j-x_i}\phi_{ij}dz\nonumber\\  &=& \left\{ \begin{array}{ll}
\wo\phi_{ij}(w)\frac{dw}{\al w}&\quad\mbox{ if }x_i>x_j\and\\
\wo\phi_{ij}(w)\frac{dw}{\al}&\quad\mbox{ if }x_i\le x_j,\\ \end{array}\right.
\eeql
with $x_1=x$ and $x_2=x'$. We take this to
define a \de{parabolic Higgs field}.  Denote the parabolic Higgs field
constructed in this way by
$\wo{\phi}$. In Simpson's language \cite{si90} is $\wo{\phi}$ just a filtered
regular Higgs field.

This defines a correspondence between Higgs $V$-bundles and parabolic Higgs
bundles (with appropriate parabolic weights).  In order to make this a
correspondence between the stable objects we simply have to check that the
invariant subbundles correspond---this is easy.  Thus we can apply many of our
preceding results to spaces of stable parabolic Higgs bundles.

\bsu{Reduction to the case $n_0 = 0$}{detred}

Suppose that at some marked points the $V$-bundle $E$ has $x = x'$ so that $n_0
> 0$.
Number the marked points so that these are the last $n_0$.  We can twist by
a line $V$-bundle to make the isotropy zero at such points.  Thus, as far as
$E$ is
concerned, the orbifold structure at these points is irrelevant and we suppose
that
$M$ only has $n-n_0$ marked points.  More precisely, we can construct $\b M$
from $M$
using the smoothing process that gives $\wo M$ but only at the last $n_0$
marked
points.  We write $\b E$ for $E$ considered as a $V$-bundle over $\b M$.

We also have to consider the canonical $V$-bundle $K$.  Notice that $K = K_{\b
M}\otimes_{i=n-n_0+1}^{n} L_{i}^{\al_i-1}$ so that there is a natural inclusion
$H^0(K^2_{\b M}) \hookrightarrow H^0(K^2_M)$ given by $ s \mapsto
s\otimes_{i=n-n_0+1}^{n} s_{i}^{2\al_i-2}$.  (Here the $L_i$ are point
$V$-bundles and $s_i$ are the canonical sections, as in \refsu{orbdiv}.)  We
identify $H^0(K^2_{\b M})$ with its image in $H^0(K^2_M)$.  From \refeq{Taylor
Higgs} it is clear that $\det(\phi)$ vanishes to order $2\al -2$ in $z$ at the
last $n_0$ marked points (since $x = x'$ there).  It follows that
$\det(\phi)\in
H^0(K^2_{\b M})$ for all Higgs fields $\phi$ on $E$.  Moreover, if we pass from
$\phi$ to $\b\phi$ by applying the smoothing process for Higgs fields at the
last $n_0$ marked points, then it is clear that $(\b E,\b\phi)$ is a Higgs
$V$-bundle over $\b M$.  Notice that by \refeq{Taylor Higgs 2} $\b\phi$ is
holomorphic at the last $n_0$ marked points because there we have $x = x'$.

The process outlined above is invertible.  For the proofs in the remainder of
this section therefore, although we will be careful to state results for $q\in
H^0(\b M,K_{\b M}^2)$ and $n_0\ge 0$, we can assume that $n_0 = 0$ without loss
of
generality.

\bsu{Generic fibres of the determinant map}{gendet}

We assume that $2g+n-n_0>3$.  Let $q \in H^0(K_{\b M}^2)$ and consider the
corresponding section $\wo q \in H^0(\wo{K_{\b M}^2})$.  We want to suppose
that
$\wo q$ has simple zeros and that none of the zeros of $\wo q$ occurs at a
marked point (of $\b M$) but first we would like to know that such behaviour is
generic.
\ble{generic} The generic section $\wo q \in H^0(\wo{K_{\b M}^2})$ has simple
zeros, none of which is at a marked point of $\b M$, provided $2g+n-n_0>3$.
\ele
\bpf
We can assume that $n_0=0$.  Notice that $\wo{K^2} = K_{\wo M}^2
\otimes_{i=1}^{n}L_{p_i}$, where $L_{p_i}=L_i^{\al_i}$ is the point bundle
associated to a marked point $p_i$.  We know that the $\wo q$ with simple zeros
form a non-empty Zariski-open set in the complete linear system $|K_{\wo M}^2
\otimes_{i=1}^{n}L_{p_i}|$.  The extra condition that none of the zeros is at a
marked point is obviously also an open condition, so we only need to check that
the resulting set is non-empty.

If $n=1$ then we only need to show that the marked point is not a base-point of
the linear system.  Similarly, if there are several marked points then it
suffices to show that none is a base point, because then the sections vanishing
at a given marked point cut out a hyperplane in the projective space $|K_{\wo
M}^2 \otimes_{i=1}^{n}L_{p_i}|$.  Using \cite[IV, proposition 3.1]{ha77}, this
is equivalent to showing that $h^0(K_{\wo M}^2 L_{p_j}^{*}
\otimes_{i=1}^{n}L_{p_i}) = h^0(K_{\wo M}^2 \otimes_{i=1}^{n}L_{p_i}) - 1$ for
each $j = 1,\dots,n$---this follows from an easy Riemann-Roch calculation,
provided $2g+n>3$.  \epf

\ble{generic2} Let $\phi$ be a Higgs field on $E$ with $\det(\phi)=q$ and $\wo
q$ generic in the
sense of \refle{generic}.  Then $\wo q$ has simple zeros at each marked point
where $x=x'$.
Moreover, at every marked point of $M$ we have $\wo\phi_{21}\ne 0$ and
$\wo\phi_{12}\ne
0$, where
$\wo\phi_{21}$ and $\wo\phi_{12}$ are as in \refeq{Taylor Higgs}.  \ele \bpf
Using
\refeq{Taylor
Higgs} we have that, in our local coordinates around a marked point,
\beql{Taylor Higgs 3} \phi =
\left( \begin{array}{rr} z^{\al-1}\wo\phi_{11}(z^\al) & z^{\al + x - x'
-1}\wo\phi_{12}(z^\al) \\
z^{x' - x -1}\wo\phi_{21}(z^\al)& -z^{\al-1}\wo\phi_{11}(z^\al)
\end{array}\right)dz,
\eeql
assuming that $x \ne x'$.  If $x' = x$ then the $(2,1)$-term is $z^{\al
-1}\wo\phi_{21}(z^\al)dz$.  Here the $\wo\phi_{ij}$ are holomorphic functions.

If $\wo q$ is generic then it is non-zero at a marked point of $\b M$ and has
at
most a simple zero at a marked point where $x=x'$---in fact there will
be a zero at such a point.  It follows that we must have that $\det(\phi) = q$
vanishes exactly to order $\al -2$ in $z$ in the first case and order $2\al -2$
in the second.  Hence $\wo\phi_{21}(0) \ne 0$ and $\wo\phi_{12}(0) \ne 0$ at
each marked point of $M$.  \epf

Henceforth we assume that $\wo q$ is a generic section, as in
\refle{generic}, and construct $\det^{-1}(q)$.  For the purposes of exposition
we also assume that $n_0 = 0$.  We face two problems in defining the
spectral variety of $\phi$ or $\wo\phi$---the first is that $\wo\phi$ has
simple
poles at the marked points and the second is that $\wo q$ is not the
determinant
of $\wo\phi$.  Let
$s_{p_i}=s_i^{\al_i}$ be the canonical section of the point-bundle $L_{p_i}$
associated to a marked point $p_i$ and let $s_0 = \otimes_{i=1}^{n}s_{p_i}$ be
the corresponding section of $\otimes_{i=1}^{n}L_{p_i}$.  Define
\beq
\oo q = \wo q s_0 \in H^0(K_{\wo M}^2
\otimes_{i=1}^{n}L_{p_i}^2) &\mbox{and}& \oo\phi = \wo\phi s_0 \in
\mbox{ParEnd}_0(\wo
E)\otimes K_{\wo
M} \otimes_{i=1}^{n}L_{p_i}.\label{eq:oophi}
\eeq
It is clear that $\det(\oo\phi)=\oo q$ and that $\oo q$ has simple zeros
(including one at each marked point).  Eventually we will need to reverse the
construction of $\oo\phi$ from $\phi$; this can be done for a given $\oo\phi
\in
\mbox{ParEnd}_0(\wo E)\otimes K_{\wo M} \otimes_{i=1}^{n}L_{p_i}$ provided
$\oo\phi$ obeys the obvious vanishing conditions at each marked point.

The square root $\sqrt{-\oo q}$ defines a smooth Riemann surface $\lift M$ with
double-covering $\pi:\lift M \to \wo M$ and branched at the zeros of $\oo q$.
Therefore there are $4g-4 + 2n$ branch-points and the Riemann-Hurwitz
formula gives the genus of $\lift M$ as $\lift g = 4g-3 + n$.  We set
$s=\sqrt{-\oo q}$---a section of $\pi^*(K_{\wo
M}\otimes_{i=1}^{n}L_{p_i})$---and $\lift\phi = \pi^*\oo\phi$.  Moreover, if
$\sigma$ is the involution interchanging the leaves of $\lift M$ then $\sigma^*
s = -s$ and $\lift\phi$ is $\sigma$-invariant.

In order to reverse the passage from $E$ to $\wo E$ we have to keep track of
the
quasi-parabolic data.  The following lemma is useful here. (Applying the
involution $\sigma$, the same result holds for $\sigma^*L =\ker(\lift\phi -
s)$.)
\ble{quasi structure}
If $\phi$ is a Higgs field on $E$ with $\det(\phi) = q$ and $\wo
q$ generic in the sense of \refle{generic}, then the kernel of $\lift\phi + s$
(with $s$,
$\lift\phi$ defined as above) is a line subbundle $L$ of $\pi^*\wo E$ and, at a
marked point (of
$\b M$) $p$, $0 \subsetneq L_{\pi^{-1}(p)} \subsetneq \pi^*\wo E_{\pi^{-1}(p)}
= \wo E_p$ describes
the quasi-parabolic structure.  \ele
\bpf
At a marked point, using \refeq{Taylor Higgs 2} and
\refeq{oophi}, we write \beql{matrix oophi} \oo\phi = \left( \begin{array}{rr}
w\wo\phi_{11}(w) &
w\wo\phi_{12}(w) \\ \wo\phi_{21}(w)& -w\wo\phi_{11}(w)
\end{array}\right)\frac{dw}\al,
\eeql with, from
\refle{generic2}, $\wo\phi_{21}(0) \ne 0$ and $\wo\phi_{12}(0) \ne 0$.
This means that $\oo\phi$ is not zero at a marked point.  Similarly, using the
fact that $\wo q$ has simple zeros, $\oo\phi$ is non-zero at {\em every}
branch
point. Now consider $\lift\phi + s$: since $\det(\lift\phi +
s)\equiv 0$ this mapping has nullity 1 or 2 at every point.  Because
$\lift\phi$
is trace-free and $s$ is scalar it follows that zeros of $\lift\phi + s$ can
only occur at zeros of $s$ \ie at the ramification points.  However, since
$\oo\phi$ is non-zero at a branch point $p$ it is impossible for $\lift\phi +
s$
to be zero at $\pi^{-1}(p)$.  So $\lift\phi + s$ is nowhere zero and the kernel
is a line bundle. Finally, if $p$ is a marked point it is clear from
\refeq{matrix oophi} that $\ker(\lift\phi + s)_{\pi^{-1}(p)}$ is spanned by
$\left( 0 , 1 \right)^T$ in our local coordinates.  The result about the
quasi-parabolic structure follows.
\epf

\bth{fibres of det} Suppose that $2g+n-n_0>3$.
Given $q \in H^0(\b M,K_{\b M}^2)$ such that $\wo q$ is generic in the sense of
\refle{generic} the
fibre of the determinant map $\det^{-1}(q)$ is biholomorphic to the Prym
variety of the covering
$\pi:\lift M \to \wo M$, determined by $q$ (via $\wo q'$).  \eth
\bpf
Since the proof is familiar \cite[theorem 8.1]{hi87} we only sketch it.  We
assume $n_0=0$.  Fix $q$ such that $\wo q$ is generic and $\lift M$ as
constructed above and also a line bundle $P$ over $\lift M$ such that
$P\sigma^* P = \pi^*(K_{\wo M}^*\La^2\wo E\otimes_{i=1}^{n}L_{p_i}^* )$.

Suppose that $(E,\phi)$ is a Higgs $V$-bundle over $M$ with $\det(\phi)=q$.
Consider the parabolic bundle $\wo E$ and $\oo\phi\in \mbox{ParEnd}_0(\wo
E)\otimes K_{\wo M} \otimes_{i=1}^{n}L_{p_i}$ with determinant $\oo q$ defined
as above.  Now set $L = \ker(\lift\phi + s)$ and notice that $L\sigma^*L \cong
\pi^*(K_{\wo M}^*\La^2\wo E\otimes_{i=1}^{n}L_{p_i}^*)$. Since $P$ was chosen
to
have the same property $LP^*$ is an element of the Prym variety.

Conversely, we consider $L$ such that $LP^*$ is a given point in the Prym
variety.  The push-forward sheaf $\pi_*\co(L)$ is locally free analytic of rank
2 and so defines a rank 2 holomorphic vector bundle $W$ over $\wo M$.  There is
a natural quasi-parabolic structure on $W^*$ at a branch point $p$ because $W_p
= (J_1L)_{\pi^{-1}(p)}$ and there is a natural filtration of jets $ 0 \subset
L^*_{\pi^{-1}(p)} \subset (J_1L)^*_{\pi^{-1}(p)}.$ The Hecke correspondence for
quasi-parabolic bundles defines a rank 2 holomorphic bundle $W'^*$:  that is,
the quasi-parabolic structure on $W^*$ defines a natural surjective map
$\co(W^*) \surjarrow \cs$, where $\cs$ is a sheaf supported at the branch
points, and the kernel of this map is locally free analytic of rank 2 and so
defines $W'^*$.

This construction of $W'$ actually recovers $\wo E$:  there is a natural map
$\co(W) \to \co(W')$ which induces an inclusion $L \hookrightarrow \pi^*W'$.
Similarly there is an inclusion $\sigma^*L \hookrightarrow \pi^*W'$.  As
subbundles of $\pi^*W'$, $L$ and $\sigma^*L$ coincide precisely on the
ramification points so that there is a map $L\oplus \sigma^*L \to \pi^*W'$
which
is an isomorphism away from the ramification points.  It follows that $\La^2W'
=
\La^2\wo E$ and that $W' = \wo E$.  Moreover, at a marked point $p$ the
inclusion $L_{\pi^{-1}(p)} \hookrightarrow \pi^*\wo E_{\pi^{-1}(p)} = \wo E_p$
gives the quasi-parabolic structure and so we recover the original
$V$-bundle $E$ (see \refpr{V-parabolic} and \refle{quasi structure}).  We
recover the Higgs field simply by defining $\lift\phi :  \pi^*\wo E \to
\pi^*(\wo E\otimes K_{\wo M}\otimes_{i=1}^{n}L_{p_i})$ by $\lift\phi(e) = \mp
se$ according as $v \in L$ or $v \in \sigma^* L$.  Since this is
$\sigma$-invariant it descends to define $\oo\phi$ on $\wo M$---this is
trace-free with determinant $\oo q$ and recovers the old $\oo\phi$.  At a
marked
point $p$, we have $\ker(\lift\phi_{\pi^{-1}(p)}) = L_{\pi^{-1}(p)}$ and
hence,
in coordinates which respect the quasi-parabolic structure, the
$(1,2)$-, $(2,2)$- and $(1,1)$-components of $\oo\phi$ vanish at $p$ to first
order in $w$.  Of course this is exactly the condition for $\oo\phi$ to define
$\wo\phi$ via \refeq{oophi} and to $\wo\phi$ there corresponds a Higgs field
$\phi$ on the $V$-bundle $E$.

Finally note that if there was an $\oo\phi$-invariant subbundle $L'$ then there
would be a section $t \in H^0(K_{\wo M}\otimes_{i=1}^{n}L_{p_i})$ such that for
any $l\in L'$, $\oo\phi(l) = tl$.  Since $\oo\phi$ is trace-free it would
follow
that $\oo q = \det(\oo\phi) = -t^2$---contradicting the assumption that $\oo q$
has simple zeros.  So $\oo\phi$ has no invariant subbundles and the same is
therefore true of $\wo\phi$ and $\phi$.  \epf
Notice that this shows that a Higgs field in the generic fibre of $\det$ leaves
no sub-$V$-bundle invariant (compare \refsu{higalg}).

\bsu{The case $g=n-n_0=1$}{detspe}

We briefly indicate how the preceding arguments can be modified to identify the
generic fibre of the determinant map when $g=n-n_0=1$.  We outline the argument
working with $V$-bundles although the proofs again require translation to the
parabolic case.  As before we simplify the exposition by supposing that $n_0=0$
so that there is a single marked point $p=p_1$.

\ble{invariants} If $g=n-n_0=1$ then every Higgs field has an
invariant sub-$V$-bundle.  \ele \bpf Since $h^0(K^2) = 1$ the natural squaring
map $H^0(K) \to H^0(K^2)$ is surjective.  Thus, given any Higgs field $\phi$,
$\det(\phi) = -s^2$ for some $s \in H^0(K)$.  Consider $\theta_\pm = \phi \pm
s$:  if $\phi\ne 0$ this is non-zero (if $\phi = 0$ then there is nothing to
prove) but has determinant zero and so we have line $V$-bundles $L_\pm
\hookrightarrow E$ with $L_\pm \subseteq \ker \theta_\pm$.  Clearly $L_\pm$ are
invariant, with $\phi$ acting on $L_\pm$ by multiplication by $\mp s$.  \epf

Since the squaring map is surjective, \refle{generic} certainly can't hold in
this case---we now consider any non-zero determinant to be `generic'.  Using
\refle{invariants} we see that any Higgs field with a generic (\ie non-zero)
determinant has two invariant sub-$V$-bundles $L_+$ and $L_-$.

Notice that $K = L_1^{\al_1 -1}$ and so sections of $K$ are multiples of the
canonical section
$s_1^{\al_1-1}$ and those of $K^2$ are multiples of $s_1^{2\al_1-2}$.  Thus in
\refeq{Taylor Higgs 3} $\wo\phi_{11}(z^{\al_1})$ and exactly one of
$\wo\phi_{12}(z^{\al_1})$ and $\wo\phi_{21}(z^{\al_1})$ are non-zero at the
marked point, while the other must vanish to first order in $w=z^{\al_1}$.  A
small local calculation using \refeq{Taylor Higgs 3} shows that $L_+$ and
$L_-$ have the same isotropy; it is $x$ if $\wo\phi_{21}(0)=0$ and $x'$ if
$\wo\phi_{12}(0) = 0$. Hence $L_+L_- \cong \Lambda L_{1}^{x - x'}$ or
$\Lambda L_1^{x' - x - \al}$, where the isotropy of $L_\pm$ is $x$ in the first
case and $x'$ in the second.  Using these and stability, we calculate that
$c_1(L_\pm) = r/2 + x/\al$ or $(r-1)/2 + x'/\al$, respectively, where $c_1(E) =
r + (x + x')/\al$.  Notice that the parity of $r$ determines the isotropy of
$L_\pm$.  Thus a point in the generic fibre gives a point not of a Prym variety
but of the Jacobian $\mbox{Jac}_0{M}\cong T^2$ corresponding to $L_+$.
Reversing the correspondence as in \refth{fibres of det} yields the following
result.
\bpr{special} If $g=n-n_0=1$ then for $q\in H^0(\b M,K^2_{\b
M})\setminus \{ 0 \}$ the fibre $\det^{-1}(q)$ is biholomorphic to the Jacobian
torus.  \epr

\bsu{Non-stable $V$-Bundles in Fibres of the Determinant Map}{detnon}

We have a natural inclusion of the cotangent bundle to the moduli of stable
$V$-bundles in to the
moduli of stable Higgs $V$-bundles and we would like to show, following
\cite[\S 8 ]{hi87}, that in
fact we
have a fibrewise compactification with respect to the determinant map.  Thus we
need to analyse the
fibres of the determinant map and check that, generically, the non-stable
$V$-bundles form
subvarieties of codimension at least 1. We wish to adapt Hitchin's argument
here
but there are additional complications and a new variant of the argument is
needed in the special case $g=n-n_0=1$.

\bpr{} Suppose that $2g+n-n_0>3$.  For fixed, generic, $q\in H^0(\b
M,K_{\b M}^2)$ let ${\rm Prym}(\lift M)$ be the Prym variety which is the fibre
of the determinant
map (\refth{fibres of det}).  Then the points of ${\rm Prym}(\lift M)$
corresponding to
non-stable $V$-bundles form a finite union of subvarieties of codimension at
least 1.
\epr
\bpf
Suppose $n_0=0$ and consider $L_E\hookrightarrow E$ a destabilising
sub-$V$-bundle, with $\wo{L_E} \hookrightarrow \wo E$ {\em parabolically}
destabilising (see \cite{fs92} and \refsu{parhig}) and $L' = \pi^*\wo{L_E}$.
The outline of the argument is similar to that of \cite[\S 8]{hi87}---with
which
we assume familiarity---but there are two problems.  Firstly, a sufficient
condition for lifts from $H^0(L'^*L^*\pi^*\La^2\wo E)$ to $H^0(L'^*\pi^*\wo E)$
to be unique is $H^0(L'^*L)=0$ but this is not always the case if $g=0$.
However, {\em invariant} lifts will still be unique because
$H^0(L'^*L)\hookrightarrow H^0(L'^*\pi^*\wo E)$ is moved by the involution
$\sigma$.  Secondly, because $\wo{L_E}$ is parabolic destabilising we can't fix
the degree of $L'^*L^*\pi^*\La^2\wo E$ in the same way that Hitchin does.  Let
the isotropy of $L_E$ be specified by an isotropy vector $(\epsilon_i)$.  A
small computation with the stability condition shows
\beq
c_1(L'^*L^*\pi^*\La^2\wo E) \le \isum\frac{\epsilon_i(x_i'-x_i)}{\al_i} + 2g
-2.
\eeq
Since $L_{\pi^{-1}(p)}$ gives the flag which describes the quasi-parabolic
structure at a marked point $p$, by \refle{quasi structure}, the subset of
$\pi^{-1}(\{p_1,\dots,p_n\})$ at which our section of $L'^*L^*\pi^*\La^2\wo E$
vanishes is just $\pi^{-1}(\{p_i :  \epsilon_i=1\})$.  Hence, for given
$(\epsilon_i)$, it is more natural to consider sections of $(\otimes_{\{i :
\epsilon_i=1\}}L_i^*)L'^*L^*\pi^*\La^2\wo E$ and these correspond to divisors
of
degree less than or equal to $\isum(\epsilon_i(x_i'-x_i)/\al_i) - n_+ +2g -2$.

For each $(\epsilon_i)$ (a finite number) we obtain a subvariety of the variety
of effective divisors and correspondingly a subvariety of the Prym variety of
codimension at least 1. \epf

\bpr{} If $g=n-n_0=1$ then
for $q\in H^0(K^2_{\b M})\setminus \{ 0 \}$ there are only a finite number of
points in the fibre
$\det^{-1}(q)$ corresponding to non-stable $V$-bundles.  \epr
\bpf
Again, we consider a destabilising sub-$V$-bundle $L_E \hookrightarrow E$ and
the corresponding parabolic bundle $\wo{L_E}$.  Since $\wo{L_E}$ is parabolic
destabilising $2c_1(\wo{L_E}) \ge c_1(\wo E) + 1$ or $2c_1(\wo{L_E}) \ge
c_1(\wo
E)$, according to whether $L_E$ has isotropy $x$ or $x'$.
Recall (from \refsu{detspe}) that $E$ has two $\phi$-invariant sub-$V$-bundles
$L_\pm$ and so is an extension $0 \to L_\pm \to E \to
L_\pm^*\Lambda \to 0$. Set $r=c_1(\wo E)$.  The discussion in
\refsu{detspe} also shows that if $r$ is even then $c_1(\wo{L_\pm})=r/2$,
$L_\pm$ have isotropy $x$ and $\wo{L_+}\wo{L_-}\cong \La^2\wo E$, while if $r$
is odd then $c_1(\wo{L_\pm}) = (r-1)/2$, $L_\pm$ have isotropy $x'$ and
$\wo{L_+}\wo{L_-}\cong \La^2\wo EL_{p_1}^*$.

Consider the sequence of bundles
\beq
0 \to \wo{L_E}^*\wo{L_\pm} \to
\wo{L_E}^*\wo E \to \wo{L_E}^*\wo{L_\pm}^*\La^2\wo E \to 0.
\eeq
and the first three terms of the associated cohomology long exact sequence.  By
assumption $H^0(\wo{L_E}^*\wo E)$ is non-zero so at least one of
$\wo{L_E}^*\wo{L_+}$ and $\wo{L_E}^*\wo{L_+}^*\La^2\wo E$ must have a non-zero
section and the same is true with $L_-$ in place of $L_+$.  If we had that
$H^0(\wo{L_E}^*\wo{L_\pm})\ne 0$ and $H^0(\wo{L_E}^*\wo{L_\pm}^*\La^2\wo E)=0$
then the inclusion of $\wo{L_E}$ in $\wo E$ would have to factor through that
of
$\wo{L_\pm}$, which is impossible as $\wo{L_\pm}$ does not destabilise.  So
$\wo{L_E}^*\wo{L_+}^*\La^2\wo E$ and $\wo{L_E}^*\wo{L_-}^*\La^2\wo E$ must have
non-zero sections.  However, considering cases according to the parity of $r$
and the isotropy of $L_E$, we see that $c_1(\wo{L_E}^*\wo{L_\pm}^*\La^2\wo
E)\le
0$.  It follows that a non-stable $V$-bundle occurs only if $\wo{L_E} \cong
\wo{L_\pm}^*\La^2\wo E$.  Since $\wo{L_+}\wo{L_-}\cong \La^2\wo E$
or $\wo{L_+}\wo{L_-}\cong \La^2\wo EL_{p_1}^*$, it follows that
$\wo{L_+}^2\cong
\La^2\wo E$ or $\wo{L_+}^2\cong \La^2\wo EL_{p_1}^*$.  Hence, if a non-stable
$V$-bundle occurs then $\wo{L_+}$ is one of the $2^{2g}=4$ possible square
roots of a given line bundle.  \epf

\bse{Representations and Higgs $V$-bundles}{rep}

Throughout this section $E\to M$ is a complex rank 2 $V$-bundle over an
orbifold
Riemann surface of negative Euler characteristic. We also suppose that a fixed
metric
and Yang-Mills connexion, $A_\Lambda$, are given on $\Lambda$.

\bsu{Stable Higgs $V$-bundles and Projectively Flat Connexions}{repsta}

Suppose that $E$ is given a Higgs $V$-bundle structure with Higgs field $\phi$,
compatible with $A_\La$.  Given a Hermitian metric on $E$ inducing the fixed
metric on $\Lambda$, there is a unique Chern connexion $A$ compatible with the
holomorphic and unitary structures and inducing $A_\Lambda$ on $\Lambda$.  The
metric also defines an adjoint of $\phi$, $\phi^*$.  Set \beq D = \partial_A +
\phi + \o\partial_A + \phi^*; \eeq this is a (non-unitary) connexion with
curvature $F_D = F_A + [\phi,\phi^*]$ and $D$ is projectively flat
if and only if the pair $(A,\phi)$ is \ymh.  The
determinant-fixing condition on $D$ is simply that it induces the fixed ({\em
unitary}) Yang-Mills connexion $A_\Lambda$ in $\Lambda$.

Conversely, given a connexion $D$ (with fixed determinant) and a Hermitian
metric on $E$, inducing the fixed metric on $\Lambda$, we can decompose $D$
into
its $(1,0)$- and $(0,1)$-parts; $D=\partial_1 + \o\partial_2$.  There are then
uniquely defined operators $\o\partial_1$ and $\partial_2$ (of types $(0,1)$
and
$(1,0)$ respectively) such that $d_1=\partial_1 + \o\partial_1$ and
$d_2=\partial_2 + \o\partial_2$ are unitary connexions.  Define $\phi =
(\partial_1-\partial_2)/2$ and $d_A=(d_1+d_2)/2$ so that $\o\partial_A =
(\o\partial_1+ \o\partial_2)/2$.  Clearly $(A,\phi)$ is a Higgs pair if and
only
if $\o\partial_A(\phi) = 0$, \ie $\phi$ is holomorphic; if we define $D''=
\o\partial_A + \phi$ then this condition becomes $D''^2 = 0$.  Here $D''$ is a
first order operator which satisfies the appropriate $\o{\partial}$-Leibniz
rule.  Moreover, if $D''^2=0$ then $(A,\phi)$ is \ymh\ if and only if $D$ has
curvature $- \pi \ci\, c_1(\Lambda)\Omega I_E$.

{}From now on suppose that $D$ has curvature $- \pi \ci\,c_1(\Lambda)\Omega
I_E$.
We call a Hermitian metric (with fixed determinant) \de{twisted harmonic} with
respect to $D$ if the resulting $D''$-operator satisfies $D''^2=0$.  Using the
fact that the curvature of $D$ is $- \pi \ci\,c_1(\Lambda)\Omega I_E$, a small
calculation shows that the condition for the metric to be twisted harmonic is
$F_1 = F_2$, where $F_i$ is the curvature of $d_i$, for $i=1,2$.  If the metric
is twisted harmonic then $D''$ defines a Higgs $V$-bundle with respect to which
the metric is \hymh.  Clearly the processes of passing from a Higgs $V$-bundle
to a projectively flat connexion and vice-versa are mutually inverse and
respect
the determinant-fixing conditions.

We prove an existence result for twisted harmonic metrics,
following \cite{do87}.  The connexion $D$ on $E$ comes from a projectively flat
connexion
in the corresponding principal $GL_2(\C)$ $V$-bundle $P$ with
$E=P\times_{GL_2(\C)}\C^2$.
Hence $D$ determines a holonomy representation $\rho_D:\pi_1^V(M) \to
PSL_2(\C)$.  Let ${\rm Herm}^+_2$ denote the $2\times 2$ positive-definite
Hermitian matrices (with the metric described in \cite[\S VI.1]{ko87}).  The
corresponding $V$-bundle of Hermitian metrics on $E$ is just
$H'=P\times_{GL_2(\C)}{\rm Herm}_2^+$.  Here $GL_2(\C)$ acts on ${\rm
Herm}^+_2$
by $h \mapsto \o{g}^T h g$, for $h\in {\rm Herm}_2^+$ and $g\in GL_2(\C)$.
This
is an action of $PSL_2(\C)$ and so $H'$ is flat and can be written as
$H'=H'_{\rho_D}={\cal H}^2 \times_{\rho_D} {\rm Herm}^+_2$ (where ${\cal H}^2$
is the universal cover of $M$).  A choice of Hermitian metric on $E$ is a
section of $H'_{\rho_D}$ or a $\pi_1^V(M)$-equivariant map ${\cal H}^2 \to
{\Herm}_2^+$---is this map is harmonic in the sense that it minimises energy
among such maps?

Using the determinant-fixing condition, we suppose that the map to ${\rm
Herm}^+_2$ has constant determinant 1. We identify the subspace of ${\rm
Herm}_2^+ \cong GL_2(\C)/U(2)$ in which the image of the map lies with
$SL_2(\C)/SU(2) \cong {\cal H}^3$.  So we consider sections of the flat ${\cal
H}^3$ $V$-bundle $H_{\rho_D}={\cal H}^2 \times_{\rho_D} {\cal H}^3$:  the
sections of $H_{\rho_D}$ are precisely the types of map considered by Donaldson
in \cite{do87}.  The condition that a metric $h$ be twisted harmonic will then
be precisely that it is given by a harmonic $\pi_1^V(M)$-equivariant map $\lift
h:  {\cal H}^2 \to {\cal H}^3$.

Donaldson shows that the Euler-Lagrange condition for the map $\lift h$ to be
harmonic is just $d_A^*(\phi + \phi^*) = 0$ and moreover that, at least in the
smooth case and when $\rho_D$ is irreducible, such a harmonic map always
exists.
This Euler-Lagrange condition agrees with our definition of a twisted harmonic
metric.  For the existence of such harmonic maps we either follow Donaldson's
proof directly or argue equivariantly, as in \refsu{ymhequ}, obtaining the
following results.
\bpr{Donaldson} Let $\rho_D:\pi_1^V(M) \to PSL_2(\C)$ be an irreducible
representation and $s_0$ a section of the flat ${\cal H}^3$ $V$-bundle
$H_{\rho_D}={\cal H}^2 \times_{\rho_D} {\cal H}^3$.  Then $H_{\rho_D}$ admits a
twisted harmonic section homotopic to $s_0$.  \epr

\bco{twisted exists} Let $\Lambda$ have a fixed Hermitian metric and compatible

Yang-Mills connexion.  Given an irreducible $GL_2(\C)$-connexion $D$ on $E$
with
curvature $- \pi \ci\, c_1(\Lambda)\Omega I_E$ and fixed determinant, $E$
admits
a Hermitian metric of fixed determinant which is twisted harmonic with respect
to $D$.  Hence $D$ determines a stable Higgs $V$-bundle structure on $E$ with
fixed determinant, for which the metric is \hymh.  \eco

\bco{} Let $E$ have a fixed Hermitian metric and let $\Lambda$ have a
compatible
Yang-Mills connexion.  Let $D$ be an irreducible $GL_2(\C)$-connexion on $E$
with curvature $- \pi \ci\, c_1(\Lambda)\Omega I_E$ and fixed determinant.
Then
there is a complex gauge transformation $g \in {\cal G}^c$, of determinant 1,
such that the {\em fixed} metric is twisted harmonic with respect to $g(D)$.
Hence $g(D)$ determines a stable Higgs $V$-bundle structure on $E$ with fixed
determinant.  \eco

To identify the space of such projectively flat connexions modulo gauge
equivalence with our moduli space of Higgs $V$-bundles we have to consider the
actions of the gauge groups and the question of irreducibility.  We have the
following result adapted from \cite[theorem 9.13 \& proposition 9.18]{hi87}.
\bpr{irreducibles} Let
$E\to M$ be a complex rank 2 $V$-bundle with a fixed Hermitian metric and
compatible Yang-Mills
connexion on the determinant line $V$-bundle $\Lambda$.  Then the following
hold.
\begin{enumerate} \item\label{irreducible'} A \ymh\ pair $(A,\phi)$ (with fixed
determinant) is
irreducible if and only if the corresponding projectively flat
$GL_2(\C)$-connexion $D=\partial_A +
\dbar_A + \phi + \phi^*$ is irreducible.  \item\label{gauge} Two irreducible
$GL_2(\C)$-connexions
on $E$ with curvature $- \pi \ci\, c_1(\Lambda)\Omega I_E$ (and fixed
determinant), $D$ and
$D'$, are equivalent under the action of ${\cal G}^c$ if and only if the
corresponding
\ymh\ pairs $(A,\phi)$ and $(A',\phi')$ are equivalent under the action of
$\cal G$.
\end{enumerate} \epr

\bsu{Projectively Flat Connexions and Representations}{reprep}

In the smooth case projectively flat connexions are described by
representations
of a universal central extension of
the fundamental group (see \cite{hi87}, also \cite[\S 6]{ab82}).  However over
an orbifold
Riemann surface there is in general no {\em one} central extension which will
do \cite[\S
3]{fs92} but the determinant-fixing condition tells us that the
appropriate central extension to use is the fundamental group of the circle
$V$-bundle
$S(\Lambda)$.  Let $(y_i)$ ($0\le y_i \le \al_i-1$) denote the
isotropy of a line $V$-bundle $L$ and let $b = c_1(L) - \isum (y_i/\al_i)$.
The orbifold fundamental group of $S(L)$ is
well-known (see, for instance, \cite[\S 2]{fs92}) and has presentation
\beql{circle}
\begin{array}{rcl} \pi_1^V(S(L)) &=& \langle a_1,b_1,\dots, a_g,b_g,q_1,\dots,
q_n,h \quad | \\ &
&\quad [a_j,h]=1, \ [b_j,h]=1,\ [q_i,h]=1,\ q_i^{\al_i} h^{y_i}=1,\ q_1\dots
q_n[a_1,b_1]\dots[a_g,b_g]h^{-b}=1 \rangle.  \end{array} \eeql

\bpr{representations} Let $\Lambda\to M$ be a line $V$-bundle with a fixed
Hermitian
metric and compatible Yang-Mills connexion.  Let $S(\Lambda)$ be the
corresponding circle
$V$-bundle.  Then there is a bijective correspondence between \begin{enumerate}
\item
conjugacy classes of irreducible representations $\pi_1^V(S(\Lambda)) \to
SL_2(\C)$ such
that
the generator $h$ in \refeq{circle} is mapped to $-I_2\in SL_2(\C)$ and \item
isomorphism
classes of pairs $(E,D)$, where $E$ is a $GL_2(\C)$ $V$-bundle with $\Lambda^2E
= \Lambda$
and $D$
is an irreducible $GL_2(\C)$ connexion on $E$ with curvature $- \pi \ci\,
c_1(\Lambda)\Omega
I_E$ and inducing the fixed connexion on $\Lambda$.  \end{enumerate} \epr
\bpf The proof
can be carried over from \cite[theorem 4.1]{fs92} (compare also
\cite[theorem 6.7]{ab82}) except that we need to replace $U(2)$ with
$GL_2(\C)$ at each stage---only the unitary structure on the determinant line
$V$-bundle is
necessary for the proof.  \epf

Since \refpr{representations} insists that $h$ maps to
$-I_2$, it is sufficient to consider a central $\Z_2$-extension rather than the
central
$\Z$-extension of $\pi_1^V(M)$ given by the presentation \refeq{circle}---this
is
equivalent to adding the relation $h^2=1$ to that presentation. Then it is only
the {\em
parity} of the integers $y_i$ and $b$ that matters. Something a little subtler
is true.
Recall \refre{roots}: it is sufficient to consider topological
$\Lambda$'s modulo the equivalence $\Lambda \sim \Lambda L^2$. Moreover, the
topology of
$\Lambda$ is specified by the $y_i$'s and $b$ (\refpr{v-bundles})---write
$\Lambda=\Lambda_{(b,(y_i))}$ to emphasise this. Clearly if $(b,(y_i)) \equiv
(b',(y'_i))
\pmod2$ (meaning that the congruence holds componentwise) then
$\Lambda_{(b,(y_i))} \sim
\Lambda_{(b',(y'_i))}$. However, if $\al_i$ is odd then $L$ can be chosen so
that
tensoring
by $L^2$ brings about a change $y_i \mapsto y_i+1$; if any $\al_i$ is even then
a change
$b
\mapsto b+1$ is possible similarly. These equivalences correspond to
group isomorphisms between the corresponding presentations \refeq{circle}, with
the added
relation $h^2=1$. Thus we normalise the $y_i$'s and $b$ to find exactly
one representative of each class, supposing that
\beql{normal}
y_i = \left\{\begin{array}{rl}
0 &\mbox{if $\al_i$ is odd;}\\
0,1&\mbox{if $\al_i$ is even;}
\end{array}\right.
\quad b = \left\{\begin{array}{rl}
0 &\mbox{if at least one $\al_i$ is even;}\\
0,1&\mbox{if no $\al_i$ is even.}
\end{array}\right.
\eeql
This is equivalent to considering only the following \de{square-free}
topological $\Lambda$'s:
\beql{normal2}
\Lambda \in \left\{\begin{array}{rl}
\{\otimes_{\al_i {\rm\ even}}L_i^{\delta_i} : \delta_i =0,1 \} &\mbox{if at
least
one $\al_i$ is even;}\\
\{L^{\delta_0} : \delta_0 = 0,1\} &\mbox{if no $\al_i$ is even,}
\end{array}\right.
\eeql
where $L$ has no isotropy with $c_1(L)=1$ and
the $L_i$ are the point $V$-bundles of \refsu{orbdiv}.

An alternative way to understand these $\Z_2$-extensions of the fundamental
group is as
follows.
Since $SL_2(\C)$ double-covers $PSL_2(\C)$ any representation
$\rho_D:\pi_1^V(M)\to PSL_2(\C)$
induces a central $\Z_2$-extension of $\pi_1^V(M)$:
\beql{extensor}
0 \to \Z_2 \to \G \to \pi_1^V(M) \to 0.
\eeql
Since the group of central $\Z_2$-extensions of $\pi_1^V(M)$ is discrete, the
$\G$ thus induced is constant over connected components of the representation
space.

So, given any $\rho_D$, we obtain an extension $\Gamma$:  what invariants
$(b,(y_i))$ characterise these $\Gamma$'s and thus the central
$\Z_2$-extensions
of $\pi_1^V(M)$?  The answer is that $(b,(y_i))$ can be supposed to have one of
the normalised forms given by \refeq{normal} and so these parameterise the
central $\Z_2$-extensions of $\pi_1^V(M)$.  This is because the image of each
generator of \refeq{circle} has exactly two possible lifts to $SL_2(\C)$ except
that $h$ must map to $-I_2$:  choosing lifts at random, the relations
$q_i^{\al_i} h^{y_i}=1$ and $q_1\dots q_n[a_1,b_1]\dots[a_g,b_g]h^{-b}=1$ of
\refeq{circle} will be satisfied for exactly one choice of normalised
$(b,(y_i))$.  By our previous discussion, this is exactly equivalent to
considering only the square-free $\Lambda$'s of \refeq{normal2}.

As well as topological types of determinant line $V$-bundles we need to
consider
topological types
of rank 2 $V$-bundles with the {\em same} determinant line
$V$-bundle---\refpr{representations}
deals with all topological types of $V$-bundles with the same determinant line
$V$-bundle
simultaneously.  These types can be determined following the ideas of \cite[\S
4]{fs92}, as
follows.  The various topological types are distinguished by the rotation
numbers associated to the
images of the elliptic generators $q_i$ of the presentation \refeq{circle}.  By
this we mean that
the image of $q_i$ has conjugacy class described by the roots of its
characteristic polynomial, necessarily of the
form $e^{\pi\ci r_i/\al_i}$, $e^{-\pi\ci r_i/\al_i}$, for $0\le r_i \le \al_i$;
these
$r_i$ are the
\de{rotation numbers}.  Notice that the relation $q_i^{\al_i} h^{y_i}=1$ means
that $r_i$
has the same parity as $y_i$ and this is the only {\em a priori}
restriction on $r_i$. Call an abstract set of rotation numbers $(r_i)$
\de{compatible
with $\Lambda$} if $r_i$ has the same parity as $y_i$. The result
is the following and the proof, using \refpr{v-bundles}, is easy.
\ble{rotation numbers} The
topological types
of $GL_2(\C)$ $V$-bundles $E$ with fixed determinant constructed in
\refpr{representations}
correspond to the rotation numbers $r_i$ associated to the images of the
elliptic
generators $q_i$
of the presentation \refeq{circle}. \ele

Denote the space of representations of $\pi_1^V(S(\Lambda))$ into $SL_2(\C)$,
sending the generator $h$ of \refeq{circle} to $-I_2$, by
$\Hom^{-1}(\pi_1^V(S(\Lambda)),SL_2(\C))$ and the {\em irreducible}
representations by $\Hom^{*,-1}(\pi_1^V(S(\Lambda)),SL_2(\C))$, for a fixed
line
$V$-bundle $\Lambda$.  For any set of rotation numbers $(r_i)$ (with $0\le r_i
\le \al_i$ and $r_i\equiv y_i \pmod2$) we have a corresponding subset
$\Hom^{-1}_{(r_i)}(\pi_1^V(S(\Lambda)),SL_2(\C))$ and, by
\refpr{representations} and the results of \refsu{repsta}, a bijection between
$\Hom^{*,-1}_{(r_i)}(\pi_1^V(S(\Lambda)),SL_2(\C))/SL_2(\C)$ and the moduli
space of stable Higgs $V$-bundles (with fixed determinants) on the topological
$E$ corresponding to the rotation numbers (\refle{rotation numbers}).

The representation space
$\Hom^{*,-1}_{(r_i)}(\pi_1^V(S(\Lambda)),SL_2(\C))/SL_2(\C)$ can be thought of
as the quotient of a set of $2g+n$ matrices subject to conditions corresponding
to the relations of \refeq{circle} and so has a natural topology; whether
this description makes it into a smooth manifold is by
no means immediate.  Therefore we use the bijection with
the moduli space of stable Higgs $V$-bundles, which is easily seen to be a
homeomorphism, to define a manifold structure on this representation space.  In
summary we have the following theorem.
\bth{} Let $M$ be an orbifold Riemann surface with negative Euler
characteristic.  Let $\Lambda$ be a fixed line $V$-bundle over $M$ and $(r_i)$
a
set of rotation numbers compatible with $\Lambda$.  Then the representation
space $\Hom^{*,-1}_{(r_i)}(\pi_1^V(S(\Lambda)),SL_2(\C))/SL_2(\C)$ is a complex
manifold of dimension $6(g-1)+2(n-n_0)$, where $n_0$ is the number of rotation
numbers congruent to 0 (mod $\al$).  \eth
\bre{twist again} In \refre{roots} we noted that twisting by a non-trivial
topological root $L$ induces a map $\cm(E,A_\Lambda) \leftrightarrow
\cm(E\otimes L,A_\Lambda)$, preserving the topology of $\Lambda$ but altering
that of $E$.  On the level of representations there is an equivalent map.
Given
any element $\lift\rho_D \in \Hom^{-1}_{(r_i)}(\pi_1^V(S(\Lambda)),SL_2(\C))$
we
can obtain a representation with different rotation numbers and covering the
same $PSL_2(\C)$-representation, by altering the signs of the images of certain
of the generators of \refeq{circle}.  We can change the sign of
$\lift\rho_D(q_i)$ (bringing about a change of rotation number $b_i \mapsto
\al_i-b_i$) provided $\al_i$ is even and provided an even number of such
changes
is made---these conditions preserve the relations $q_i^{\al_i} h^{y_i}=1$ and
$q_1\dots q_n[a_1,b_1]\dots[a_g,b_g]h^{-b}=1$.  \ere

When there are no reducible points we can apply, among other results,
\refpr{metric} and \refco{topology}. By \refle{rotation numbers} we can
discuss
the existence of reducible points in
terms of the rotation numbers.  (Either $\Lambda$ or a specific set of
rotation numbers may provide an obstruction to the existence of reductions.)
The discussion in \refsu{ymhmod} shows that reductions exist if and
only if there exists an isotropy vector $(\epsilon_i)$ such that
\beq \isum\frac{\epsilon_i(x'_i-x_i) + (x'_i +
x_i)}{\al_i}\equiv c_1(\Lambda) \pmod 2. \eeq
A small calculation expresses this in terms of the rotation numbers. Thus we
obtain the following result.
\bpr{} Let $M$ be an orbifold Riemann
surface with negative Euler characteristic.  Let $\Lambda$ be a fixed line
$V$-bundle over
$M$ with isotropy $(y_i)$ and $c_1(\Lambda) = b + \isum(y_i/\al_i)$. Let
$(r_i)$ be a
compatible set of rotation numbers. Then the representation space
$\Hom^{-1}_{(r_i)}(\pi_1^V(S(\Lambda)),SL_2(\C))/SL_2(\C)$ contains reducible
points if
and only if there exists an isotropy vector
$(\epsilon_i)$ such that
\beq
\isum\frac{\epsilon_ir_i}{\al_i}\equiv b \pmod 2.
\eeq
When no reducible points exist the complex manifold
$\Hom^{-1}_{(r_i)}(\pi_1^V(S(\Lambda)), SL_2(\C))/SL_2(\C)$
\begin{enumerate}
\item admits a complete
hyper-K\"ahler metric and
\item is connected and simply-connected.
\end{enumerate}
\epr

\bsu{Real Representations}{reprea}

In the previous subsection we discussed $SL_2(\C)$-representations of central
extensions of the orbifold fundamental group.  Here we study the submanifold of
$SL_2(\R)$-representations.  First notice that any irreducible representation
into $SL_2(\C)$ can fix at most one disk ${\cal H}^2\subset {\cal H}^3$ because
the intersection of two fixed disks would give a fixed line and hence define a
reduction of the representation.  Moreover, any representation which does fix a
disk can be conjugated to a real representation and the conjugation action of
$SL_2(\C)$ then reduces to that of $SL_2(\R)$.

Now consider the action of complex conjugation on a representation. Recall
that,
via \refpr{representations} and \refco{twisted exists}, irreducible
representations
correspond to stable Higgs $V$-bundles. Note that $\pi_1^V(S(\Lambda))$ and
$\pi_1^V(S(\o\Lambda))$ are isomorphic via the map $h \mapsto h^{-1}$: the
following proposition follows, exactly as in \cite{si90}.
\bpr{} Let $E$ be a complex rank 2 $V$-bundle such that $\Lambda$ has a fixed
Hermitian
metric and compatible Yang-Mills connexion.  Let
$\lift\rho_D:\pi_1^V(S(\Lambda)) \to
SL_2(\C)$ be an irreducible
representation, sending $h$ to $-I_2$, with corresponding stable Higgs
$V$-bundle structure on $E$, $(E_A,\phi)$. Then the complex conjugate
representation
(thought of as a representation of $\pi_1^V(S(\o\Lambda))$) determines a Higgs
$V$-bundle
structure on $\o E$, isomorphic to $(E_A,-\phi)^*$.  \epr

\bco{} Let $E$ be a complex rank 2 $V$-bundle such that $\Lambda$ has a fixed
Hermitian
metric and
compatible Yang-Mills connexion.  Let $\lift\rho_D$ be an irreducible {\em
real} representation
$\lift\rho_D:\pi_1^V(S(\Lambda)) \to SL_2(\R)$, sending $h$ to $-I_2$, with
corresponding Higgs $V$-bundle structure $(E_A,\phi)$.  Then there
is an isomorphism of Higgs $V$-bundles $(E_A,\phi)\cong (E_A,-\phi)$.  \eco

Consider the involution on the moduli space of stable Higgs $V$-bundles (with
fixed unitary structure and determinants) defined by $\sigma:  (E,\phi) \mapsto
(E,-\phi)$, where now $E$ denotes a holomorphic $V$-bundle and $(E,\phi)$ is a
stable Higgs $V$-bundle.  The fixed points of $\sigma$ can be determined much
as
the fixed points of the circle action were in the proof of \refth{Morse}.  If
$(E,\phi)$ is itself fixed then $\phi=0$ and $E$ is a stable $V$-bundle.
Suppose now that $\phi\ne 0$.  If $(E,\phi)$ is only fixed up to complex
gauge-equivalence then we have an element $g \in {\cal G}^c$ such that
$g(E,\phi) = (E,-\phi)$.  Since $g$ fixes $E$ it must fix the Chern connexion
$A$ and since $g$ cannot be a scalar it leads to a reduction of $A$ to a direct
sum of $U(1)$-connexions.  Hence we have a holomorphic decomposition $E = L
\oplus L^*\La$, where, without loss of generality, we may suppose that $2c_1(L)
- c_1(\La) \ge 0$.  Since $(A,\phi)$ is an irreducible pair, $g$ must have
order
2 in ${\cal G}^c$ and fix $A$.  It follows that with respect to this
decomposition (or, if $A$ has stabiliser $SU(2)$, {\em choosing} a
decomposition) we can write
\beq g =\pm\left( \begin{array}{ll} \ci & 0\\ 0 &
-\ci \end{array} \right) &\mbox{and}& \phi = \phantom{\pm}\left(
\begin{array}{rr} t & u\\ v & -t \end{array} \right).  \eeq
(Since our Higgs $V$-bundle is stable, we must have $v$ non-zero.)  Calculating
the conjugation-action of $g$ on $\phi$ we find that $t=0$.

Recall that we chose $L$ with $2c_1(L) - c_1(\La) \ge 0$ but to avoid
semi-stable points
(when $u=0$) we suppose that there is strict inequality.  Exactly as in the
proof of
\refth{Morse} we consider the topological possibilities
$L=L_{(m,(\epsilon_i))}$:
we can have any $(m,(\epsilon_i))$ such that $2c_1(L) > c_1(\Lambda)$ and
$c_1(\wo{KL^{-2}\La}) = r \ge 0$.  Then the possible holomorphic structures and
the values
of $v$, modulo the $\C^*$ automorphism group, are given by the effective
(integral)
divisors of order $r$ and taking square roots.  A difference from
\refth{Morse}
is that $u$ needn't be zero; indeed, $u$ can take any value in
$H^0(KL^2\La^*)$.  We
obtain the following result, where $l$ is defined as in \refsu{orbint}.
\bpr{real
manifolds}
Let $M$ be an orbifold Riemann surface of negative Euler characteristic and
suppose that
$E\to M$ admits no reducible \ymh\ pairs. Then the fixed points of the
involution
induced on $\cm(E,A_\La)$ by the mapping $(A,\phi) \mapsto (A,-\phi)$ consist
of complex
$(3g-3+n-n_0)$-dimensional submanifolds $\cm_0$ and $\cm_{(m,(\epsilon_i))}$,
for every
integer $m$ and isotropy vector $(\epsilon_i)$ such that
\beq
l < 2m + \isum \frac{\epsilon_i(x'_i-x_i)}{\al_i} \le l + 2g -2 +
\isum\frac{\epsilon_i(x_i'-x_i)}{\al_i} + n_-.
\eeq
The manifold $\cm_0$ is the moduli space of stable $V$-bundles with fixed
determinants, while $\cm_{(m,(\epsilon_i))}$
is a rank $(2m -l + g -1  +n_+)$ vector-bundle over a
$2^{2g}$-fold covering of $S^r\wo
M$, where $r = l -2m + 2g -2  + n_-$.
\epr

We interpret this as a result about $PSL_2(\R)$-representations of
$\pi_1^V(M)$. Again, a representation $\rho_D$ of $\pi_1^V(M)$ into $PSL_2(\R)$
induces a
central $\Z_2$-extension $\Gamma$ of $\pi_1^V(M)$, as in \refeq{extensor},
which is just $\pi_1^V(S(\Lambda))$ with the added relation $h^2=1$, for some
square-free
$\Lambda$. Consider the points of $\Hom^{-1}(\pi_1^V(S(\Lambda)),SL_2(\R))$
covering $\rho_D$. On the level of
representations there are $2^{2g+n_2-1}$ (or $2^{2g}$ if $n_2 = 0$) choices of
sign for
the images of certain generators and these correspond to twisting a stable
Higgs
$V$-bundle by any of the $2^{2g+n_2-1}$ (or $2^{2g}$) holomorphic roots of the
trivial
line $V$-bundle. In particular, if $n_2 \ge 1$ then the topology of the
associated $E$ is
only determined up to twisting by the $2^{n_2-1}$ non-trivial topological roots
(see
\refre{twist again}).

Excluding the topologically non-trivial roots, we have an action of
$\Z_2^{2g}$ on the fixed point submanifolds of \refpr{real manifolds} which is
easily
seen to be free if $E$ admits no reducible \ymh\ pairs. Moreover, even when
$E$
admits reducibles there will be fixed submanifolds $\cm_{(m,(\epsilon_i))}$
with
\beq
l \le 2m + \isum \frac{\epsilon_i(x'_i-x_i)}{\al_i} \le l + 2g -2 +
\isum\frac{\epsilon_i(x_i'-x_i)}{\al_i} + n_-,
\eeq
exactly as in \refpr{real manifolds}, and the actions of $\Z_2^{2g}$ on these
will be free provided the first inequality is strict.

The quantity $2m -l + \isum \{\epsilon_i(x'_i-x_i)/\al_i\}  =
2c_1(L_{(m,(\epsilon_i))}) -
c_1(\Lambda)$ is just the Euler class of the flat $\R\P^1$ $V$-bundle
$S(\rho_D) = S(L_{(m,(\epsilon_i))}^2 \Lambda^*)$ associated to the
$PSL_2(\R)$-representation (this is well-defined as it is invariant under
twisting $E$
by non-trivial topological roots). Note that, just as it is possible to have
topologically
distinct line $V$-bundles with the same Chern class, it is possible to have
topologically
distinct $\R\P^1$ $V$-bundles with the same Euler class---they are
distinguished by their
isotropy. The central $\Z$-extensions of $\pi_1^V(M)$
induced by the universal covering $\widetilde{PSL_2\R} \to PSL_2\R$ are just
the (orbifold) fundamental groups of the flat $\R\P^1$ $V$-bundles
$S(\rho_D)$ (see \cite{jn85}). Using the above discussion and the method of
\refpr{real manifolds}, we obtain the following result (compare \cite{jn85})
and, as a corollary, a Milnor-Wood inequality.
\bpr{psl2r reps}
Let $M$ be an orbifold Riemann surface of negative Euler characteristic. For
$\rho_D$ a $PSL_2(\R))$-representation of $\pi_1^V(M)$ let
$\Hom_{\rho_D}(\pi_1^V(M),PSL_2(\R))$ denote
the corresponding connected component. Let $(y_i)$ be the isotropy and $b
+ \isum (y_i/\al_i)$ the Euler class of the associated flat $\R\P^1$ $V$-bundle
$S(\rho_D)$. Provided $b + \isum (y_i/\al_i) > 0$,
$\Hom_{\rho_D}(\pi_1^V(M),PSL_2(\R))/PSL_2(\R)$ is a smooth complex
$(3g-3+n-n_0)$-dimensional manifold, diffeomorphic to a rank $(g - 1 + b + n
-n_0)$
vector-bundle over $S^{2g-2-b}\wo M$.
\epr

\bco{milnor-wood}
Let $M$ be an orbifold Riemann surface of negative Euler characteristic. Then
the Euler
class $b + \isum (y_i/\al_i)$ of any flat $PSL_2(\R)$ $V$-bundle
satisfies
\beq
|b + \isum \frac{y_i}{\al_i}| \le 2g -2 + n -\isum\frac1{\al_i}.
\eeq
\eco
\bpf
In \refpr{psl2r reps} we must have $b \le 2g-2$. The result follows since $y_i
\le
\al_i-1$.
\epf

\bsu{Teichm\"uller Space for Orbifold Riemann Surfaces}{reptei}

Assume, as usual, that $M$ is an orbifold Riemann surface of negative Euler
characteristic.  For a Fuchsian group such as $\pi_1^V(M)$, Teichm\"uller
space,
denoted $\ct(M)$, is the space of faithful representations onto a discrete
subgroup of $PSL_2\R$ modulo conjugation (see Bers's survey article
\cite{be72}).  Our previous results allow us to identify Teichm\"uller space
with a submanifold of the moduli space.

Let $\ct_{-4}(M)$ denote the space of orbifold Riemannian metrics of constant
sectional curvature -4, modulo the action of the group of diffeomorphisms
homotopic to the identity, $\cd_0(M)$.  There is a bijection between
$\ct_{-4}$ and $\ct$ as each metric of constant negative curvature determines
an
isometry between the universal cover of $M$ and $\ch^2$ and hence a faithful
representation of $\pi_1^V(M)$ onto a discrete subgroup of $PSL_2\R$ and
conversely each such representation realises $M$ as a geometric quotient of
$\ch^2$.

The results of \cite{jn85}, as well as those of \cite[\S 11]{hi87}, suggest
that
Teichm\"uller space is the component of the real representation space taking
the
extreme value in the Milnor-Wood inequality, \refco{milnor-wood}. Working with
the
holomorphic description, the results of the
previous subsection show that the extreme is achieved when $E=L\oplus
L^*\Lambda$ with
$L^2\Lambda^*$ having the topology of $K$ and a holomorphic structure such that
$\wo{KL^{-2}\Lambda}$ has sections: in other words we must have $L^2\Lambda^* =
K$
(holomorphically). We suppose then that $E=K \oplus 1$ ($\Lambda^2E$ can be
normalised to be square-free but this is not necessary).
The corresponding Higgs field is just \beq \phi = \left( \begin{array}{cc} 0 &
u\\ v & 0\\ \end{array}\right), \eeq where $u\in H^0(K^2)$ and $v\in
\C\setminus
\{0\}$.  There is a $\C^*$-group of automorphisms so that we can normalise with
$v=1$.

Exactly as in \cite[theorem 11.2]{hi87}, we can identify Teichm\"uller
space with the choices of $u$ \ie with $H^0(K^2)$.  The two preliminaries which
we need are the strong maximum principle for orbifolds (the proof is entirely
local and generalises immediately; see \cite{jt80}) and the following orbifold
version of a theorem of Sampson \cite{ee69}.
\bpr{Sampson}
Given two orbifold Riemannian metrics of constant
sectional curvature -4 on $M$, $h$ and $h'$, there is a unique element of
$\cd_0$ which is a harmonic map between $(M,h)$ and $(M,h')$.  \epr
\bpf This is a reformulation of
\refpr{Donaldson}.  The metrics $h$ and $h'$ give two discrete, faithful
representations of
$\pi_1^V(M)$ into $PSL_2\R$, one of which we consider fixed and the other we
denote $\rho'$.  The identity map on $M$ lifts to an
orientation-preserving diffeomorphism $g$
of $\ch^2$ which is equivariant with respect to the actions of the two
representations.
Taking this $g$ as an initial section of the $V$-bundle $H_{\rho_D}={\cal H}^2
\times_{\rho'} {\cal
H}^3$ of \refpr{Donaldson} (via the inclusion $\ch^2 \subset \ch^3$) we obtain
a harmonic section
$g'$ homotopic to $g$.  This is real and defines a harmonic diffeomorphism
between $(M,h)$
and $(M,h')$. As $g'$ is homotopic to $g$ the resulting harmonic
diffeomorphism
is homotopic to the identity.
Uniqueness follows either by a direct argument or from uniqueness over $\lift
M$, where
$\lift M$ is as in \refco{smooth covering}.  \epf

We obtain the following theorem, which agrees with classical results due to
Bers and others \cite{be72}.
\bth{ball} Let $M$ be an orbifold Riemann surface of negative Euler
characteristic.  Let
$\ct(M)$ be the
Teichm\"uller space of the Fuchsian group $\pi_1^V(M)$ and $\ct_{-4}(M)$ the
space of
orbifold Riemannian
metrics on $M$ of constant sectional curvature -4, modulo the action of the
group of
diffeomorphisms homotopic to the identity.
Then $\ct(M)$ and
$\ct_{-4}(M)$ are homeomorphic to $H^0(K^2)$, the space of holomorphic
(orbifold)
quadratic differentials on $M$. Hence Teichm\"uller space is homeomorphic to
$\C^{3g-3+n}$. \eth

We conclude by considering orbifold Riemannian metrics in greater detail.
Considered as a metric on the underlying Riemann surface, $\wo
M$, an orbifold Riemannian metric $h$ on $M$ has `conical singularities' at
the marked points. To see this recall that locally
$M$ is like $D^2/\Z_\al$ with
$h$ a $\Z_\al$-equivariant metric on $D^2$. If $c_h(r)$ denotes the
circumference of a
geodesic circle of radius $r$ about the origin in $D^2$ (with respect to $h$),
then
$lim_{r\to
0}(c_h(r)/r)=2\pi$. Since this circle covers a circle in $D^2/\Z_\al$ exactly
$\al$ times
the metric on the quotient has a \de{conical singularity} at the origin, with
\de{cone
angle} $2\pi/\al$.

Consider a Riemannian metric on $M$ which, near a marked point $D^2/\Z_\al$, is
compatible with the complex structure and so has the form $h(z) dz \otimes d\o
z$. If we
set $w=z^\al$, then $w$ is a local holomorphic coordinate on $\wo M$.   We find
that
the resulting `Riemannian metric' on $\wo M$ is given by \beq
\frac{h(w^{1/\al})}{\al^2|w|^{2(1-1/\al)}} dw \otimes d\o w. \eeq Notice that
$h(w^{1/\al})$
is well-defined by the $\Z_\al$-equivariance of $h$.  This `Riemannian metric'
has a singularity like
$|w|^{-2(1-1/\al)}$ at the origin and is compatible with the complex structure
away from
there.  Hence we obtain a compatible \lq singular Riemannian metric' on $\wo
M$:  the
induced metric on $\wo M$ is continuous and induces the standard topology.

How does such a singular Riemannian metric compares with a
(smooth) Riemannian metric on $\wo M$?  Suppose that $g$ is a fixed
Riemannian metric on $\wo M$,
compatible with the complex structure.  Since $\wo M$ is compact any two
Riemannian
metrics give metrics on $\wo M$ which are
mutually bounded and so will be equivalent for our purposes---we may as well
use the Euclidean
metric in any local chart.  Now, $h$ and $g$ will give mutually bounded metrics
on any compact
subset of $M\setminus \{ p_1,\dots,p_n \}$.  However, for small Euclidean
distance $r$ from $p$, the
singular metric has distance like $r^{1/\al}$. These are
exactly the types of
singularities of metrics considered by McOwen and Hulin-Troyanov in
\cite{mo88,ht92}: they
consider metrics which satisfy $h/g = O(r_g^{2k})$ as $r_g(z) = d_g(0,z) \to
0$,
for some $k\in (-1,\infty)$.  As McOwen points out, our \lq singular Riemannian
metrics' have exactly this form with $k =-1 + 1/\al$.  Interpreting
\refco{negative curvature} in the light of this discussion we obtain the
following result.  (Our result is weaker than McOwen's since we consider only
$k
=-1 + 1/\al$ but the case of general $k\in (-1,\infty)$ can be obtained by a
limiting argument as in \cite{ns93}.)
\bthn{McOwen, Hulin-Troyanov}{conical} Let $\wo M$ be a Riemann surface with
marked points $\{p_1,\dots,p_n\}$ with orders of isotropy
$\{\al_1,\dots,\al_n\}$.  If the
genus $g$ and orders of
isotropy satisfy \beq 2-2g-n+\sum_{i=1}^n 1/\alpha_i < 0 \eeq then $\wo
M\setminus
\{p_1,\dots,p_n\}$ admits a unique compatible Riemannian metric $h$ of constant
sectional curvature
-4 such that, for $i=1,\dots,n$, $h$ has a conical singularity at $p_i$ with
cone angle
$2\pi/\al_i$.  \eth

\end{document}